\documentclass[trsc,nonblindrev]{informs3} 

\OneAndAHalfSpacedXI 

\usepackage{dirtytalk}
\usepackage{multirow}
\usepackage[official]{eurosym}
\usepackage{soul}
\usepackage{xltabular}
\usepackage[ruled]{algorithm2e}
\usepackage[table,xcdraw]{xcolor}

\usepackage{natbib}
 \bibpunct[, ]{(}{)}{,}{a}{}{,}%
 \def\bibfont{\small}%
\TheoremsNumberedThrough     

\EquationsNumberedThrough    

\begin{document}


\RUNAUTHOR{Nicolet and Atasoy}

\RUNTITLE{A Choice-Driven Service Network Design and Pricing Including Heterogeneous Behaviors}

\TITLE{A Choice-Driven Service Network Design and Pricing Including Heterogeneous Behaviors}

\ARTICLEAUTHORS{%
\AUTHOR{Adrien Nicolet, Bilge Atasoy}
\AFF{Dept. of Maritime and Transport Technology, Delft University of Technology, The Netherlands, \\ \EMAIL{\{A.Nicolet, B.Atasoy\}@tudelft.nl}}
} 

\ABSTRACT{%
The design and pricing of services are two of the most important decisions faced by any intermodal transport operator. The key success factor lies in the ability of meeting the needs of the shippers. Therefore, making full use of the available information about the demand helps to come up with good design and pricing decisions. With this in mind, we propose a Choice-Driven approach, incorporating advanced choice models directly into a Service Network Design and Pricing problem. We evaluate this approach considering three different mode choice models: one deterministic with 4 attributes (cost, time, frequency and accessibility); and two stochastic also accounting for unobserved attributes and shippers' heterogeneity respectively. To reduce the computational time for the stochastic instances, we propose a predetermination heuristic. These models are compared to a benchmark, where shippers are solely cost-minimizers. Results show that the operator's profits can be significantly improved, even with the deterministic version. The two stochastic versions further increase the realized profits, but considering heterogeneity allows a better estimation of the demand.
}%


\KEYWORDS{Network Design, Pricing, Mode Choice, Heterogeneity}
\pagebreak
\maketitle

\section{Introduction}\label{intro}

In intermodal freight transport, planning at the tactical level is of key importance to make the best use of existing infrastructure and available assets and to ensure reliable transport plans. An appropriate way of managing this task is through Service Network Design (SND) problems, as they cover most of the tactical decisions~\citep{crainic2000service}. They can support the decisions of intermodal operators about the itineraries to be served, the offered frequencies and how demand should be assigned to these services.

Until recently, pricing was not explicitly covered in most SND models although it plays a crucial role in the success of the planning~\citep{tawfik2018pricing,li2015pricing}. As pointed out by~\citet{macharis2004opportunities}, intermodal transport pricing is a difficult task as costs must be accurately computed and some knowledge of the market situation has to be gained. Indeed, the costs faced by an intermodal carrier are various~\citep{li2005medium}: some of them, e.g. crew costs or contracts with infrastructure manager, are perfectly known by the carrier but other variable costs are set by external companies, such as terminal operators for the handling costs or energy suppliers for the fuel costs. For the latter, not only do they depend on external actors, but also on the transport demand as they increase together with the carried load. Although carriers have some control on the quantity of transported freight (via contract binding, for example), demand remains mostly stochastic in nature~\citep{combes2013shipment}. As a result, variable costs can only be estimated from the expected transport demand.

Regarding the pricing decision itself, some knowledge about the targeted demand, such as the willingness to pay or the transport requirements, is also of key importance. Indeed, the cost of transportation is among the main drivers of shippers' mode choice. It would, however, be inadequate to consider that shippers are purely \say{cost-minimizers} as other factors (e.g., transport time, offered quality, service frequency) play a role in the decision process, see for example~\citet{arencibia2015modelling} or~\citet{ben2013discrete}. On top of that, these factors and their importance can vary from shipper to shipper and the final decision of choosing a mode also depends on the available alternatives, hence making the planning and pricing process even more complex. On the other hand, there exists a great variety of mode choice models (see~\citet{de2014mode} for a comprehensive review) that can be used to support the planning of intermodal carriers. For example,~\citet{duan2019freight} include values of time and reliability, that are estimated from a stated preference survey, within the cost minimization of a SND model. This represents a step towards the integration of shippers' preferences within the planning process.

\subsection{Illustrative example}\label{example}

\begin{figure}[!tbp]
\centering
\includegraphics[width=0.8\textwidth]{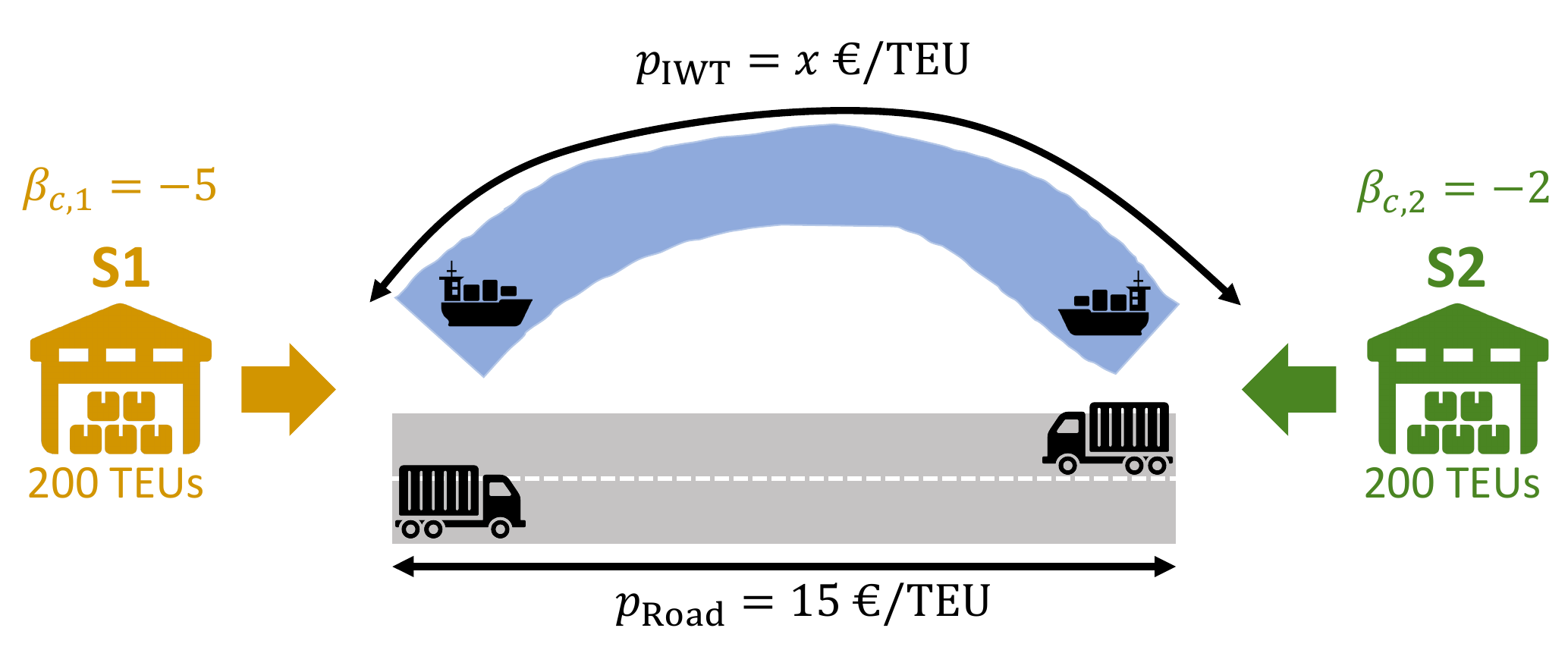}
\caption{Illustrative example with two shippers and two available transport modes.\label{fig:ToyExample}}
\end{figure}

To illustrate the benefits of using a mode choice model for the pricing decision, we consider the case in Figure~\ref{fig:ToyExample}, where two shippers, S1 and S2, want to send 200 Twenty-foot Equivalent Units (TEUs) each. To do so, they have two alternatives: Road and Inland Waterway Transport (IWT). Each mode has the following utility function for each shipper $i$:
\begin{flalign*}
\begin{cases}
    V_i^{\mathrm{IWT}} & = \beta_{f}f + \beta_{c,i}p_{\mathrm{IWT}} = 1\times5 + \beta_{c,i} \times x,\\
    V_i^{\mathrm{Road}} & = \alpha_{\mathrm{Road}} + \beta_{c,i}p_{\mathrm{Road}} = 15 + \beta_{c,i} \times 15,
\end{cases}
\end{flalign*}

where $\alpha_{\mathrm{Road}}$ is the Alternative Specific Constant (ASC) for Road, equal to 15, and the ASC for IWT is normalized to 0. $p_{\mathrm{Road}}$ is the cost of the Road alternative, set to 15~\euro{}/TEU, and $\beta_{c,i}$ represents the cost sensitivity of each shipper $i$: we assume that it is -5 for S1 and -2 for S2. $\beta_{f}$ is the weight associated to the frequency of IWT services $f$, and assumed to be 1 for both shippers.

In this example, the decision-maker is the IWT carrier that wants to set up a transport service running each working day (hence: $f=5$) and to optimize their price $x$. The carrier faces a fixed cost, $c_{\mathrm{fix}}$, of 100~\euro{} per round trip and a variable cost, $c_{\mathrm{var}}$, of 1~\euro{}/TEU. Assuming that the transport demand of shippers is split according to a logit model, the carrier aims at setting a unique price so as to maximize their profits, expressed as:
\begin{flalign*}
\Pi(x) =  \sum_{i}(200 \times \frac{e^{V_i^{\mathrm{IWT}}}}{e^{V_i^{\mathrm{IWT}}}+e^{V_i^{\mathrm{Road}}}})(x-c_{\mathrm{var}}) - f \times c_{\mathrm{fix}} = \sum_{i}(200 \times \frac{e^{V_i^{\mathrm{IWT}}}}{e^{V_i^{\mathrm{IWT}}}+e^{V_i^{\mathrm{Road}}}})(x-1) - 500
\end{flalign*}

The carrier does not necessarily know the full utility specification of the shippers. Therefore, it can opt for various strategies, here we consider three of them:
\begin{enumerate}\renewcommand{\labelenumi}{\Alph{enumi})}
\item Assume that shippers are homogeneous and purely cost-minimizers, the considered utilities may then be: $V_i^{\mathrm{IWT}} = -1x$ and $V_i^{\mathrm{Road}} = -1 \times 15$ $\forall i$;
\item Make more market study to come up with the same utility functions as above, but consider that shippers are homogeneous with a mean cost sensitivity, thus: $\beta_{c,i} = -3.5$ $\forall i$;
\item Consider also the heterogeneity regarding the cost sensitivity, thus: $\beta_{c,1} = -5$ and $\beta_{c,2} = -2$.
\end{enumerate}

\begin{figure}[!h]
\centering
\includegraphics[width=0.8\textwidth]{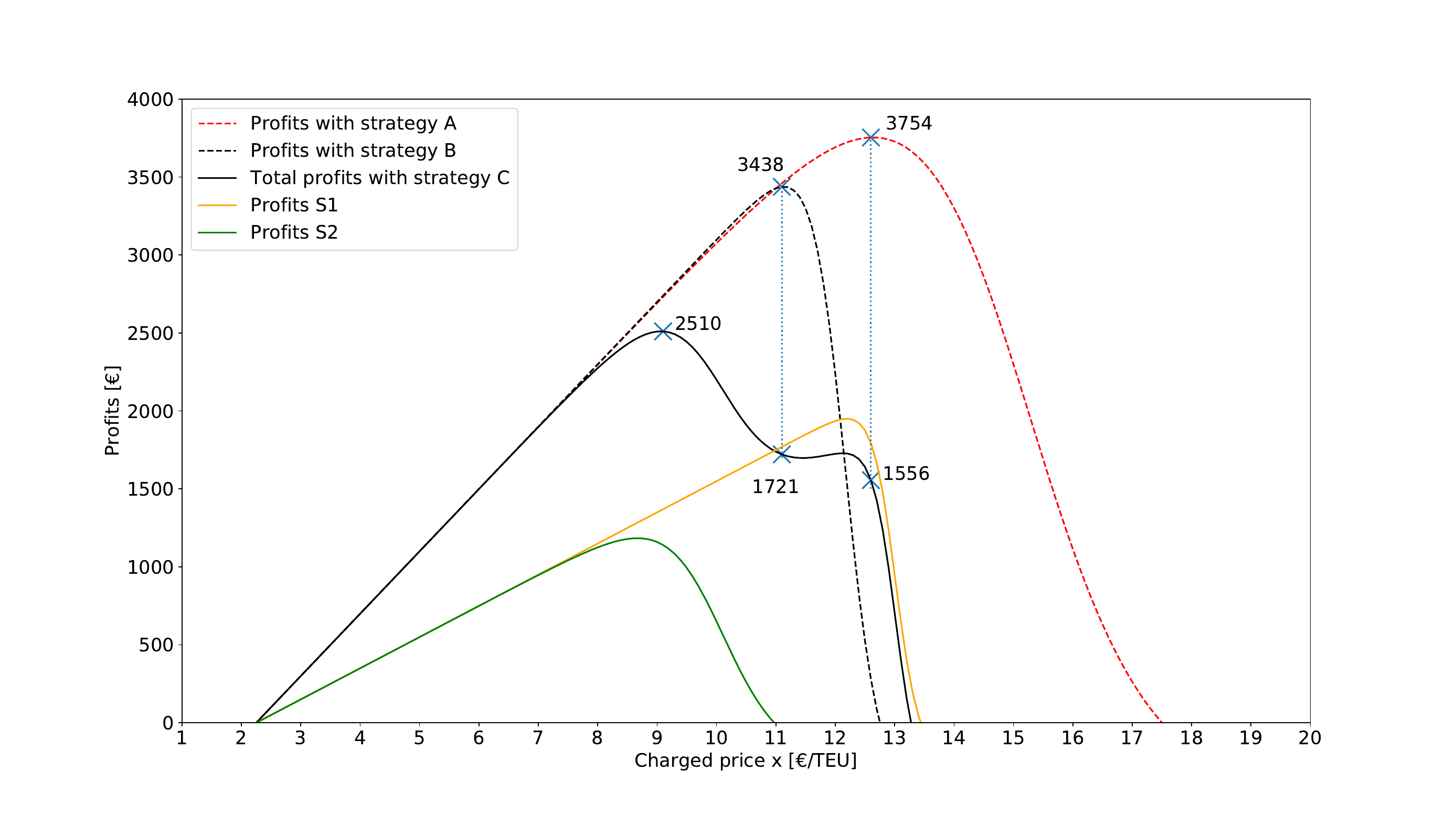}
\caption{Resulting profits for strategies A (pure cost minimization), B (homogeneous shippers), C (heterogeneous shippers) and individual profits for S1 and S2.\label{fig:ToyExample_uncap}}
\end{figure}

Assuming that the carrier has enough capacity to accommodate the demand, the resulting profits $\Pi(x)$ for each strategy are depicted in Figure~\ref{fig:ToyExample_uncap}, together with the profits stemming from each individual shipper. Before the price reaches 7~\euro{}/TEU, there is no difference between the 3 strategies: the profits grow linearly with the price between 2.25 and 7~\euro{}/TEU. This is because the price of IWT is low that its utility is much higher than Road in all three strategies: thus, the whole demand is assigned to IWT and the profits only depend on the price. The same linear relationship is observed for the individual shippers with a slope reduced by 2, since each of them represents half of the demand.

Above 7~\euro{}/TEU, the utility of IWT gradually approaches the one of Road and the relationship between price and profit becomes non-linear as the demand is also varying with the price. The highest expected profits are reached by strategy A with a price of 12.5~\euro{}/TEU. Indeed, when only considering costs, the IWT carrier can charge a higher amount as the cost of Road is relatively high and the other advantages of Road (such as low travel time) are not considered. When the price approaches and exceeds the one of Road (15~\euro{}/TEU), the profits decrease quickly as expected.

The maximal expected profits of strategy B are also particularly high and occur at 11~\euro{}/TEU. This is because the cost sensitivity of S2 is overestimated: in strategy B, it is as if the cost of Road was still too high for S2 despite the other advantages of this transport mode. But in reality, the profits stemming from S2 reach zero for a price of 11~\euro{}/TEU because Road advantages (included in $\alpha_{\mathrm{Road}}$) overcome the higher cost, thus making Road much more attractive. If strategies A or B are used for pricing, it would result in profits reduced by half when the obtained price is applied to the heterogeneous shippers.

Applying strategy C gives the highest profits, as it considers the true cost sensitivity of S2. In fact, the optimal price of 9~\euro{}/TEU almost exactly corresponds to the one maximizing the profits generated by S2 only. For higher prices, the S2 profits reduce sharply while profits from S1 only grow linearly: hence resulting in a decrease in the total profits. That is why applying the higher prices found by strategies A and B results in lower profits. It should be noted that, for strategy C, it is assumed that the carrier knows the true cost sensitivity of both shippers. In reality, it will not be the case. Nevertheless, this example showcases that the more information is known by the carrier, the more beneficial it is for the planning and pricing decisions.

A revenue management strategy would be trivial in this example with only two shippers and simple utility functions, then the optimal solution would be to set different prices for S1 and S2. However, segmentation may be difficult to identify when much more shippers are considered and less detailed information are available. In the remainder of this work, we will not consider revenue management, although we recognize that it can be an effective tool to optimize pricing decisions.

\begin{figure}[!b]
\centering
\includegraphics[width=0.8\textwidth]{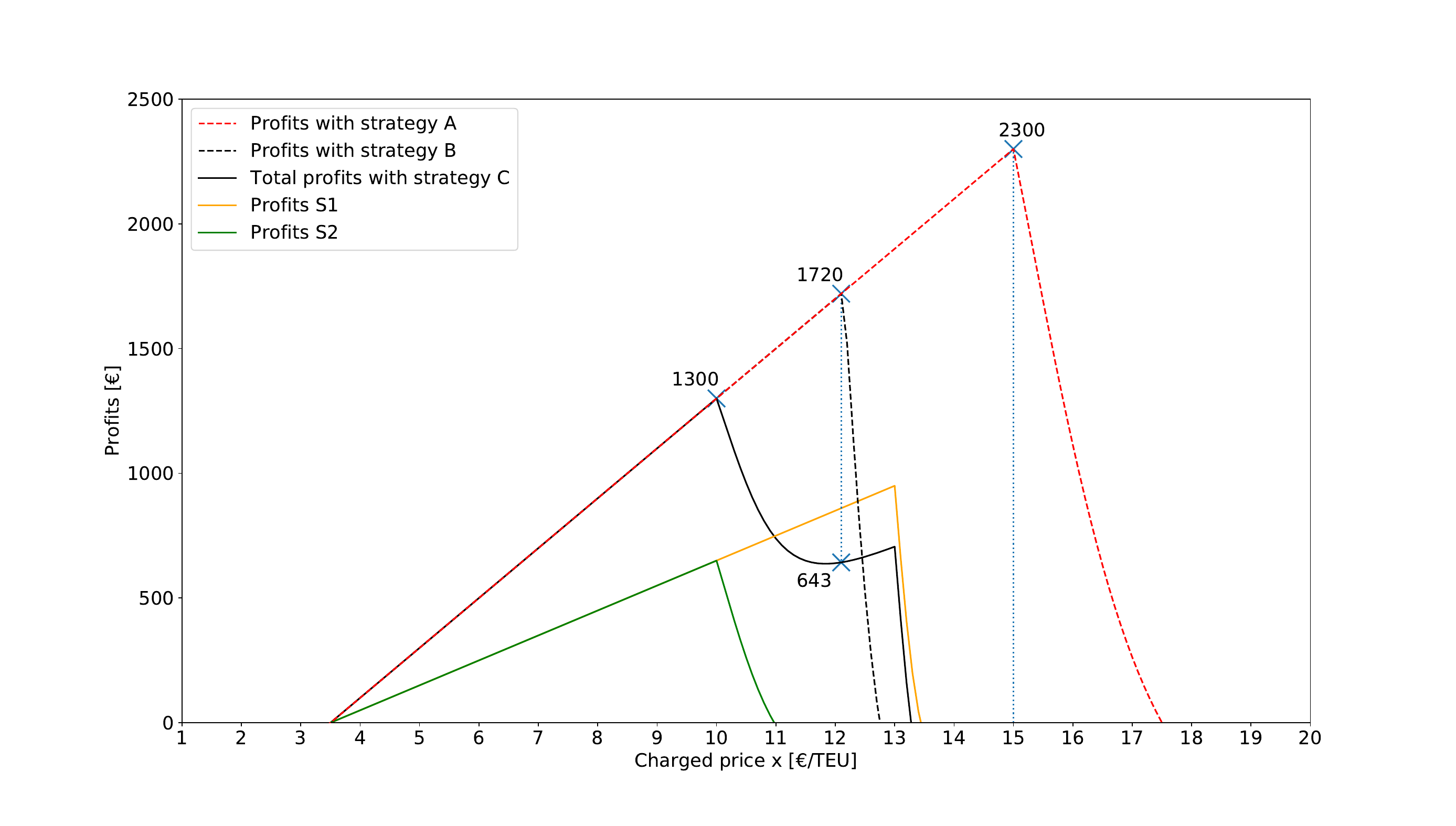}
\caption{Resulting profits for strategies A (pure cost minimization), B (homogeneous shippers), C (heterogeneous shippers) and individual profits for S1 and S2, with a vessel capacity of 20 TEUs. \label{fig:ToyExample_20TEUsCAP}}
\end{figure}

Let us conclude this example by enforcing capacity constraints for the IWT carrier. Figure~\ref{fig:ToyExample_20TEUsCAP} shows the same price-profit relationships but considering a fixed vessel capacity of 20 TEUs. Since the IWT carrier cannot serve the whole demand, the price from which positive profits are obtained is higher (3.5~\euro{}/TEU) than in the uncapacitated case. Then, the same linear relationship between price and profits is observed, but it remains valid for higher prices because the limiting factor is capacity (only 100 TEUs per direction), and not the total demand (200 TEUs per direction).

The conclusions regarding the optimal prices under different strategies remain essentially the same. However, in this capacitated case, using strategy A would result in losses for the carrier even though it has the highest expected profits. Again, the optimal price obtained through strategy C corresponds to the optimal price for S2.  

The purpose of this illustrative example is to highlight the advantages of using accurate choice models for the pricing decision: as it can substantially improve the obtained profits compared to less detailed demand representations. In fact, even if the exact parameters are not known, it is beneficial for the carrier to incorporate more information about their clients. Moreover, the impact of considering capacity constraints is demonstrated: although the previous conclusion still holds, the obtained profits are significantly reduced.

\subsection{Paper outline}\label{outline}

In this work, we aim to exploit the advantage of using advanced choice models in the planning process of an intermodal carrier. For this, we make use of \say{choice-based optimization} to combine a SND and a pricing problem with a detailed representation of shippers. Therefore, we develop a Choice-Driven Service Network Design and Pricing (CD-SNDP) model, which includes an existing mode choice model to consider shippers' behavior directly in the decision-making of the carrier.


In the rest of this paper, we review the existing literature on SNDP in Section~\ref{review}. Then, the proposed methodology is described in Section~\ref{method}, where deterministic and stochastic formulations are covered. In Section~\ref{case}, the CD-SNDP is applied to a case study and several variations of the model are compared with each other. Finally, we conclude the paper in Section~\ref{conclu} and share some insights for future research.

\section{Literature Review}\label{review}

Although it has not been applied to SNDP models yet, choice-based optimization is already used for other types of problems. Therefore, we first review the state of literature on SNDP in intermodal transport, then investigate the existing choice-driven methods in related transportation fields and finally present the main contributions of the present work.

\subsection{Service Network Design and Pricing problems in intermodal transport}\label{lit-SNDP}

The majority of existing studies on SND are formulated as a cost minimization of the transport operator and do not include the revenues of fulfilling the transport orders~\citep{elbert2020tactical,wieberneit2008service}. Nevertheless, two models using cost minimization have addressed the pricing decision. \citet{li2015pricing} determine the price charged by an intermodal carrier using a pre-defined profit margin, expressed as a given percentage of the operational costs. The price is the addition of the costs and the margin and cannot exceed a given market price. \citet{dandotiya2011optimal} include a target for the minimal profit (per transported unit) to be achieved by an intermodal operator: this translates into a constraint assuring that the applied rate is greater or equal to that target added to the operating costs. The authors also include a cost sensitivity factor representing the willingness to pay for intermodal transport rather than road and enforce that the rate difference between road and intermodal transport has to be greater or equal to this factor.

For the works applying a profit maximization, some of them do not include the pricing decision but rather assume fixed tariffs that are included as parameters into the model. \citet{andersen2009designing} apply a SND model to explore new rail services along a Polish freight corridor. The demand is represented as contracts generating a given revenue when served. The operator then decides to serve or not the contracts in order to maximize their profit. It also decides on the services' frequency and the vehicles and demand assignment to these services under vehicle balancing and capacity constraints. \citet{braekers2013optimal} are interested in designing a barge transport service along a Belgian canal considering empty container repositioning. Their SND model decides at which inland ports to stop and in which sequence as well as the fulfillment of transport demand from different clients. \citet{bilegan2022scheduled} also apply a SND model to barge transport with detailed fleet management and revenue management considerations. Different customer segments are considered as well as two different service levels (standard or express) with a given fare. The operator then decides which services to operate, what percentage of the demand to serve and how to assign the vessels and demand to the services so as to maximize their profits. The model has been developed further to include the possibility of bundling services and penalties for early and late distribution~\citep{taherkhani2022tactical}. \citet{teypaz2010decomposition} treat similar models and propose decomposition algorithms for computational efficiency. \citet{zetina2019profit} capture demand elasticity using a gravity model, where the demand is considered inversely proportional to the transport costs faced by the carrier. The decisions are whether or not an arc (or a path) is used and in which sequence to visit the demand nodes. Finally,~\citet{scherr2022stochastic} use SND to conceive a new platooning service of autonomous vehicles. They come up with a two-stage stochastic model considering scenarios to represent the demand variation. The first stage designs the services performed by \say{manually operated vehicles} and assigns rates to the different customers over all scenarios, whereas the second sets the flow of autonomous vehicles for each particular scenario.

Other works include demand functions in the profit maximization to capture the influence of prices on the transport volumes. \cite{li2005medium} design a railroad network using a concave inverse demand function. In this case, the demand for each service and each itinerary are the decision variables and the corresponding prices are computed using the inverse demand function. \citet{mozafari2011pricing} represent two competitive road carriers within a non-cooperative game model. Each carrier has to set their price so as to maximize their own profit and the demand is represented as a linear function of the carrier's price and the competitor's price. \citet{shah2012price} also investigate competition between carriers: each of them fix their price, frequency and capacity. The demand of shippers for a given carrier is represented as a function of price and frequency. The inconvenience of demand functions is that they become hard to obtain when the number of shippers or alternatives increase~\citep{li2005medium}, thus requiring a numerical estimation or some simplifying assumptions.

An increasingly common way to model SNDP problems is using Stackelberg game or bilevel programming. This formulation was first proposed in intermodal freight transport by~\citet{tsai1994optimal}. The intermodal carrier is the leader and sets the price of their services to maximize their profit. Truck carriers are followers that will adjust their prices based on the leader's decision and the exogenous demand is split between the carriers using a logit model, where the considered attributes are the prices, travel times and reliability. A general formulation for the Joint Design and Pricing (JDP) on a network has been proposed by~\citet{brotcorne2008joint}. The network operator decides on the network design and prices so as to maximize their profits. The network and rates of the competitors are assumed known and exogenous. The followers are the network users that seek to minimize their cost by selecting the services of the operator or those of the competitors. The authors propose an iterative procedure to solve the JDP. \citet{crevier2012integrated} propose a similar formulation, with the addition of capacity constraints and revenue management considerations. \citet{ypsilantis2013joint} extend the JDP formulation to include time constraints, as well as capacity constraints. Their model is used to design and price the hinterland barge services of an extended gate operator. In their work,~\citet{tawfik2019bilevel} include some level-of-service attributes in the JDP formulation. In particular, the lower level costs are more detailed as they not only consider transport costs but also the cost of capital. An iterative heuristic is later proposed to solve large instances of the JDP~\citep{tawfik2022iterative}. A similar formulation is adopted by~\citet{zhang2019integrated} to design and price rail container transport. The lower level objective is to minimize the generalized costs, made of price, transport time, convenience and security. Only the price is endogenous to the model. The same authors also propose a time varying model~\citep{zhang2019optimal,li2020reducing}. A single-level formulation is used and the demand follows a logit model with price as single attribute. The model proposed by~\citet{wang2023integrated} extends the JDP of~\citet{tawfik2019bilevel} with the introduction of additional cost components. The carrier faces some waiting costs and penalty for an under-utilisation of transport capacity, while the lower level costs also embed heterogeneous shipper classes through different values of time and reliability. 

Finally, there also exist a few different versions of Stackelberg game. A monopoly setting is proposed by~\citet{qiu2021pricing} where a hinterland carrier sets services and prices in multiple planning horizons. The followers are represented by a set of captive consignees that minimize their transport and storage costs. \citet{lee2014bi} consider three different actors as leaders and all shippers as followers. The upper level itself is represented as a three-level program where ocean carriers are leaders of terminal operators which, in turn, are leaders of land carriers~\citep{lee2014freight}. At the lower level, shippers set their production, consumption and transportation demand using \say{spatial price equilibrium}.

The relevance of bilevel models is questioned by~\citet{martin2021integrated}, especially because of the simplifying assumptions regarding demand modeling (pure cost minimizers and homogeneous preferences). They propose a SNDP model applied to an express shipping service by airplanes and trucks. In their profit maximization problem, the carrier has to set prices for some given service times that can be selected by their customers. The service time chosen by each customer is the one providing a welfare greater or equal to all the other options.

With our CD-SNDP, we show that it is possible to include advanced demand modeling, also considering heterogeneity, within a bilevel setting. In particular, the proposed formulation is inspired by the work of~\citet{tawfik2019bilevel}, where the cost minimization of shippers is replaced by the maximization of their utility. In our work, beside the costs, the utility functions also consider the transport time, the accessibility of a mode and the frequency of intermodal services. This last element implies that now, both the price and frequency decisions of the carrier have an influence on the shippers. This CD-SNDP formulation then allows for a more detailed and realistic representation of the shippers' characteristics and behavior towards the prices and services designed by the carrier. To build our model, we make use of choice-based optimization: hereafter are presented some applications of this method to other transportation problems.

\subsection{Choice-based optimization in transportation}\label{lit-CD}

The term~\say{choice-based optimization} refers to optimization problems that explicitly include a discrete choice model into their formulation~\citep{pacheco2020general}. That is why works decoupling the optimization from the demand, using iterative procedures such as simulation-optimization, are not considered here (e.g.,~\citet{liu2019framework}).

Although not for freight, choice-based optimization has been used in a few works to model passenger SND problems. \citet{wang2008multi} propose a profit maximization problem to support the design of ferry services, where the operator decides on the itineraries and schedules of the ferries. They assume that the passenger demand is split according to a logit model including two attributes: a given price, and the travel time, which is dependent on the decision variables. \citet{huang2018solving} also include a logit model into a profit maximization problem to design a car-sharing network. Among other things, the operator decides on the number of car-sharing stations to open. The utility function of car-sharing is composed of given rental costs and walking access costs. The latter are directly dependent on the number of opened stations. A drawback of these two models is that they are non-linear due to the exponential terms inherent to the logit model. A Mixed-Integer Linear Programming (MILP) including a logit mode choice model is proposed by~\citet{hartleb2021modeling} to design passenger rail services. The main decision is the selection of lines to open. To get rid of the exponential terms of the logit model, the authors precompute the modal shares of rail for each possible solution. This precomputation technique is useful when only binary or integer variables are included in the choice model. However, as mentioned by the authors, the model can become intractable when the instance size increases.


Choice-based optimization has also been applied to facility location and pricing problems. It is used by~\citet{luer2013competitive} to set up hubs and prices for an airline company. The demand is split between companies using a logit model with price as unique attribute. A similar modeling approach is adopted by~\citet{zhang2015designing} to locate retail stores and set selling prices. \citet{zhang2018game} study an intermodal dry port location and pricing problem where the route choice of shippers is determined using a logit model including six attributes, where only transport cost depends on the decision variables. The common point of these three models is that they are all non-linear: therefore, heuristics are required to solve them.

In most of the aforementioned models, the inclusion of discrete choice into an optimization problem either results in a non-linear model. In their work,~\citet{paneque2021integrating} propose a general framework to deal with more advanced choice models. In particular, the authors rely on the Sample Average Approximation (SAA) principle to deal with the non-linearities of the choice model and, therefore, come up with a MILP model. The proposed model is then applied to the pricing of parking services using a Mixed logit to represent the demand. The latter comprises price as endogenous attribute and other exogenous attributes. \citet{bortolomiol2021simulation} develop this framework further to model oligopolistic competition, whereas~\citet{schlicher2022stable} present a non-linear cooperative game to model collaborative pricing of urban mobility.

The present CD-SNDP is inspired by the work of~\citet{paneque2021integrating} to integrate an existing Mixed logit model within a bilevel setting. Specifically, error terms are included in the utilities to account for the attributes that are not captured by the model but still play a role in the mode choice. Moreover, the coefficient representing cost sensitivity is considered randomly distributed to consider the heterogeneous preferences of shippers. It is assumed that the probability distributions of the error terms and the cost coefficient are known and the resulting CD-SNDP problem is solved using stochastic optimization. The addition of these more detailed behavioral attributes within SND models aim at providing a more realistic representation of shippers' reaction to the proposed services, ultimately helping intermodal carriers to make more informed design and pricing decisions.

The following section provides a recap of the main characteristics of the previously reviewed bilevel SNDP models and sums up the contributions of our work.

\subsection{Contributions}\label{lit-contrib}

The existing bilevel models for SNDP presented in Section~\ref{lit-SNDP} are sorted in Table~\ref{tab:lit}. In particular, it shows if the models consider some competition or not and if some constraints regarding the fleet are included (e.g., size limit). It also indicates how the transport demand is modeled: most works assume that shippers are cost minimizers, while one uses price equilibrium and the last one considers demand as exogenous. Still regarding demand, it can be noticed that no existing model considers unobserved attributes that play a role in the choices of shippers. In addition, only three studies embed shippers' heterogeneity through distinct values of time (or reliability). Finally, only one work considers that frequencies also influence the demand alongside the prices endogenously in the optimization model. 

The proposed CD-SNDP is a generalization of the model by~\citet{tawfik2019bilevel}. Firstly, it generalizes the network structure as cycles and services with multiple stops are now allowed. Secondly, the shippers' objective is also generalized as they do not only proceed to a minimization of their costs, but instead maximize their utilities. These utilities contain other attributes beside the costs, such as frequency, accessibility, etc. Thirdly, our formulation generalizes the representation of shippers as it can accommodate some unobserved attributes (via randomly distributed error terms) and shippers' heterogeneity (through the Mixed logit formulation). Finally, the service frequency is made endogenous to the optimization model along with the price. A summary of the aforementioned contributions can be found on the last row of Table~\ref{tab:lit}.

\begin{table}[!tbp]
\centering
\caption{Summary of existing bilevel models for intermodal Service Network Design and Pricing problems}
\label{tab:lit}
\resizebox{\textwidth}{!}{%
\begin{tabular}{lccccccc}
\hline
\textbf{Reference} &
  \multicolumn{1}{c}{\textbf{\begin{tabular}[c]{@{}c@{}}Competition\\ included\end{tabular}}} &
  \multicolumn{1}{c}{\textbf{\begin{tabular}[c]{@{}c@{}}Fleet\\ constraints\end{tabular}}} &
  \multicolumn{1}{c}{\textbf{\begin{tabular}[c]{@{}c@{}}Demand\\ modeling\end{tabular}}} &
  \multicolumn{1}{c}{\textbf{\begin{tabular}[c]{@{}c@{}}Unobserved\\ attributes\end{tabular}}} &
  \textbf{\begin{tabular}[c]{@{}l@{}}Hetero-\\ geneity\end{tabular}} &
  \multicolumn{1}{c}{\textbf{\begin{tabular}[c]{@{}c@{}}Cross-level\\ variables\end{tabular}}} \\
\citet{tsai1994optimal} & \checkmark & \checkmark & Exogenous &  &  & Price \\
\citet{brotcorne2008joint} & \checkmark &  & Cost min. &  &  & Price \\
\citet{crevier2012integrated} & \checkmark & \checkmark & Cost min. &  & \checkmark & Price \\
\citet{ypsilantis2013joint} & \checkmark & \checkmark & Cost min. &  &  & Price \\
\citet{lee2014bi} & \checkmark & \checkmark & Price eq. &  & \checkmark & Price \\
\citet{tawfik2019bilevel} & \checkmark &  & Cost min. &  &  & Price \\
\citet{zhang2019integrated} & \checkmark & \checkmark & Cost min. &  &  & Price \\
\citet{qiu2021pricing} &  & \checkmark & Cost min. &  &  & Price \\
\citet{wang2023integrated} & \checkmark &  & Cost min. &  & \checkmark & Price \& Freq. \\
Proposed CD-SNDP  & \checkmark & \checkmark & Mixed logit & \checkmark & \checkmark & Price \& Freq. \\
\hline
\end{tabular}%
}
\end{table}

\section{Methodology}\label{method}

As previously mentioned, the present work is inspired by the bi-level JDP formulation proposed by~\citet{tawfik2019bilevel}. In both the JDP and the proposed model, the upper level represents the decisions of the transport operator and the lower level corresponds to the shippers. Just like in the JDP, the upper level consists of determining the frequency and price of services that maximize the operator's profit. However, the lower level now represents the utility maximization of the shippers, whereas in the JDP, it is assumed that shippers proceed to a minimization of their logistics costs. This change of paradigm brings additional complexity to the problem as the two decision variables of the upper level are now included in the lower level, unlike the JDP where only the price is included but not the frequency. Moreover, the competition of the transport operator is now represented as more than one alternative.

Concerning the upper level, it differs from the JDP in two aspects. Firstly, the path-based formulation of services is replaced by a cycle-based formulation. The latter is deemed more accurate to represent realistic decision-making. Indeed, most intermodal transport services go back and forth on an itinerary with a defined schedule. The cycle-based representation also enables a more elaborate representation of services as multiple intermediary stops can be added in both directions. In addition, it simplifies the asset management of the operators. In a path-based formulation, they may need to re-balance the vehicles at the end of the planning horizon; whereas a cycle-based representation ensures that each vehicle ends up at its starting point. It is noteworthy that the arc-based pricing representation is not changed compare to the benchmark. Indeed, shippers will pay only for the transport of their cargo from its origin to its destination, and not for the whole cycle.
The second difference is the addition of fleet size and cycle time feasibility constraints. The former restricts the actions of the transport operator as they do not have an infinite number of vehicles at their disposal to satisfy the demand. Moreover, the illustrative example of Section~\ref{example} showed that it significantly reduces the obtained profits. The latter determines, for each service, the number of cycles that can be performed by one vehicle during the planning horizon given the cycle's duration. A natural consequence of these additional constraints is that an heterogeneous fleet can be considered. 

In the remainder of this paper, the JDP with fleet constraints will be used as \textit{Benchmark}. The benchmark with cycle-based formulation (instead of path-based) will be further referred to as \textit{SNDP}. Finally, the proposed choice-driven model, which considers utility maximization of shippers, is denoted \textit{CD-SNDP}. The notations for the CD-SNDP are described in the following paragraphs.

\subsection{Problem formulation}\label{formul}

The transport network is represented as a directed graph $\mathcal{G} = (\mathcal{N},\mathcal{A})$, where $\mathcal{N}$ is the set of terminals and $\mathcal{A}=\{(i,j): i,j \in \mathcal{N}, i \neq j\}$ the set of links between these terminals.

\subsubsection{Upper level}\label{formul_up}

The operator's fleet is heterogeneous, therefore the different vehicle types are denoted by set $\mathcal{K}$. The number of available vehicles per type is $V_k$ and their capacity is $Q_k$.

Set $\mathcal{S}$ includes all transport services that can be run by the operator. Unlike the benchmark, where each service corresponds to a single arc of $\mathcal{A}$, a service is composed of a sequence of arcs. Each arc in this sequence is called a leg and the whole sequence of legs for a given service $s$ is denoted $\mathcal{L}_s$. The cycle-based formulation of the problem implies that the sequence starts and ends at the same node.

The maximum number of cycles of service $s$ that can be performed by vehicle type $k$ is named $W_{sk}$: it typically consists of the maximum operating time divided by the cycle time (sum of travel time and time at terminals). Each service $s$ has a fixed cost $c^{\mathrm{FIX}}_{sk}$ of operating it with vehicle type $k$ and a variable cost $c^{\mathrm{VAR}}_{ijsk}$ per container transported between terminals $i$ and $j$. Moreover, we introduce the parameter $\delta_{ijl_{s}}$, which equals one if a container traveling from $i$ to $j$ uses the service leg $l_s$ and zero otherwise.

\pagebreak

The transport operator has three decision variables in the upper level problem:
\begin{itemize}
    \item $v_{sk}$ is the number of vehicles of type $k$ that the operator allocates to each service $s$;
    \item $f_{sk}$ is the frequency of each service $s$ per vessel type $k$;
    \item $p_{ij}$ is the price per container charged to shippers requesting to transport goods from $i$ to $j$.
\end{itemize}

\subsubsection{Lower level}\label{formul_low}

The shippers are represented as a whole, i.e. their demand is aggregated. The container transport demand between terminals $i$ and $j$ is denoted $D_{ij}$. Shippers decide to assign demand to the transport operator or their competitors by the maximization of their utility. The utility function of using the services proposed by the transport operator between $i$ and $j$ is denoted $U^O_{ij}$ and is dependent on $p_{ij}$ and $f_{sk}$, whereas the utility of using a competing alternative $h$ is given as $U^h_{ij}$. 
Finally, the decision variables of the lower level consist in the number of containers that are assigned to the operator's services ($x_{ijsk}$) and to every competing alternative ($z^h_{ij}$).

All the aforementioned sets, parameters and decision variables are listed in Table~\ref{notation}.

\begin{xltabular}{\textwidth}{l  X}
\caption{Notation}\label{notation}\\
\hline \multicolumn{2}{l}{\up\textbf{Sets:}}\\

  $\mathcal{N}$ & Set of terminals (indices: $i$, $j$)\\
  $\mathcal{A}$ & Set of arcs $(i,j)$\\
  $\mathcal{K}$ & Set of vehicle types (index: $k$)\\
  $\mathcal{S}$ & Set of potential services (index: $s$)\\
  $\mathcal{L}_s$ & Set of legs of service $s \in \mathcal{S}$ (index: $l_s$)\\
  \down $\mathcal{H}$ & Set of competing alternatives (index: $h$)\\
  
 \hline
 \multicolumn {2}{l}{\up\textbf{Parameters:}}\\
 
 $V_k$ & Number of vehicles of type $k$ in the operator's fleet\\
 $Q_k$ & Capacity of vehicle type $k$ [TEUs]\\
 $W_{sk}$ & Maximum number of cycles of service $s$ that can be performed by vehicle type $k$\\
 $c^{\mathrm{FIX}}_{sk}$ & Fixed cost of operating service $s$ with vehicle type $k$ [\euro{}]\\
 $c^{\mathrm{VAR}}_{ijsk}$ & Variable cost of transport between $i$ and $j$ with service $s$ and vehicle type $k$ [\euro{}/TEU]\\
 $\delta_{ijl_{s}}$ & Dummy param. equal to 1 if container traveling from $i$ to $j$ uses service leg $l_s$, 0 otherwise\\
 $D_{ij}$ & Aggregated transport demand of shippers between $i$ and $j$ [TEUs]\\
 $U^O_{ij}$ & Utility of using the operator's services between $i$ and $j$\\
 \down $U^h_{ij}$ & Utility of using competing alternative $h$ between $i$ and $j$\\
 
 \hline
 \multicolumn {2}{l}{\up\textbf{Variables:}}\\

 $v_{sk}$ & Number of vehicles of type $k$ assigned to service $s$ by the operator\\
 $f_{sk}$ & Frequency of service $s$ operated with vehicle type $k$\\
 $p_{ij}$ & Price charged by the operator to shippers wanting to transport goods from $i$ to $j$ [\euro{}/TEU]\\
 $x_{ijsk}$ & Cargo volume using service $s$ operated with vehicle type $k$ between $i$ and $j$ [TEUs]\\
 \down $z^h_{ij}$ & Cargo volume using competing alternative $h$ between $i$ and $j$ [TEUs]\\

 \hline
\end{xltabular}

\subsubsection{Mathematical model}\label{bilevel_model}

The proposed SNDP is expressed as a Bi-Level Mixed Integer Problem (BLP):
\begin{flalign}
& \mathbf{(BLP)} \max_{v,f,p,x,z} \quad \sum_{(i,j) \in \mathcal{A}}\sum_{s \in \mathcal{S}}\sum_{k \in \mathcal{K}} p_{ij}x_{ijsk} - \sum_{s \in \mathcal{S}}\sum_{k \in \mathcal{K}} c^{\mathrm{FIX}}_{sk}f_{sk} - \sum_{(i,j) \in \mathcal{A}}\sum_{s \in \mathcal{S}}\sum_{k \in \mathcal{K}} c^{\mathrm{VAR}}_{ijsk}x_{ijsk} && \label{eq:obj_up}
\end{flalign}
\begin{flalign}
& \mathrm{s.t.} \quad \sum_{s \in \mathcal{S}} v_{sk} \enspace \leq \enspace V_k && \forall k \in \mathcal{K} && \label{eq:fleet_con} \\
& f_{sk} \enspace \le \enspace W_{sk}v_{sk} && \forall s \in \mathcal{S}, \enspace \forall k \in \mathcal{K} && \label{eq:cycle_con} \\
& \sum_{(i,j) \in \mathcal{A}} \delta_{ijl_{s}}x_{ijsk} \enspace \le \enspace Q_k f_{sk} && \forall l_s \in \mathcal{L}_s, \enspace \forall s \in \mathcal{S}, \enspace \forall k \in \mathcal{K} && \label{eq:cap_con} \\
& x_{ijsk} \enspace \le \enspace \sum_{l_s \in \mathcal{L}_s} \delta_{ijl_{s}}D_{ij} && \forall (i,j) \in \mathcal{A}, \enspace \forall s \in \mathcal{S}, \enspace \forall k \in \mathcal{K} && \label{eq:serv_con}\\
& p_{ij} \ge 0 && \forall (i,j) \in \mathcal{A} && \label{eq:p_con}\\
& v_{sk} , f_{sk} \in \mathbb{N} && \forall s \in \mathcal{S}, \enspace \forall k \in \mathcal{K} && \label{eq:vf_con}
\end{flalign}
where $x$ and $z$ solve:
\begin{flalign}
& \max_{x,z} \quad \sum_{(i,j) \in \mathcal{A}} \left( \sum_{s \in \mathcal{S}}\sum_{k \in \mathcal{K}} U^O_{ij}x_{ijsk} + \sum_{h \in \mathcal{H}} U^h_{ij}z^h_{ij} \right) && \label{eq:obj_low}
\end{flalign}
\begin{flalign}
& \mathrm{s.t.} \quad \sum_{s \in \mathcal{S}}\sum_{k \in \mathcal{K}} x_{ijsk} + \sum_{h \in \mathcal{H}} z^h_{ij} \enspace = \enspace D_{ij} && \forall (i,j) \in \mathcal{A} && \label{eq:demand_con}\\
& x_{ijsk} \ge 0 && \forall (i,j) \in \mathcal{A}, \enspace \forall s \in \mathcal{S}, \enspace \forall k \in \mathcal{K} && \label{eq:x_con}\\
& z^h_{ij} \ge 0 && \forall (i,j) \in \mathcal{A}, \enspace \forall h \in \mathcal{H} && \label{eq:z_con}
\end{flalign}

At the upper level, the objective function of the transport operator~\eqref{eq:obj_up} is to maximize their profit. It is computed as the revenues from the transported containers minus the fixed and variable costs of the offered services. Constraint~\eqref{eq:fleet_con} is the fleet size constraint for each vehicle type. Constraint~\eqref{eq:cycle_con} ensures that the service's frequency is inferior to the maximum number of cycles that can be performed by the assigned vehicles. Constraint~\eqref{eq:cap_con} assures that the total number of containers transported on each leg of every service does not exceed the available capacity of the service, whereas constraint~\eqref{eq:serv_con} ensures that no container can be assigned to a service that does not go through the origin or destination terminal of the container. The domains of the operator's decision variables are defined by constraints~\eqref{eq:p_con}-\eqref{eq:vf_con}.

Regarding the lower level, shippers seek to maximize their utility~\eqref{eq:obj_low} by assigning their containers either to the operator's services or to the competition. Moreover, constraint~\eqref{eq:demand_con} enforces the total transport to be met. Finally, constraints~\eqref{eq:x_con}-\eqref{eq:z_con} define the domain of the decision variables of the shippers.

\subsection{Model transformation}\label{linear}

The proposed bi-level problem can be reformulated as a single level problem and then linearized, for more details on these procedures the reader is referred to~\citet{tawfik2019bilevel}. For the reformulation, additional variables $\lambda_{ij}, \enspace \forall (i,j) \in \mathcal{A}$ are introduced: they represent the dual variables related to constraints~\eqref{eq:demand_con}. The model can then be transformed using the Karush-Kuhn-Tucker conditions. After this process, the following constraints appear:
\begin{flalign}
& \lambda_{ij} \enspace \le \enspace - U^O_{ij} && \forall (i,j) \in \mathcal{A}, \enspace \forall h \in \mathcal{H} && \label{eq:lU_con}\\
& \lambda_{ij} \enspace \le \enspace - U^h_{ij} && \forall (i,j) \in \mathcal{A}, \enspace \forall h \in \mathcal{H} && \label{eq:lUh_con}\\
& (-U^O_{ij} - \lambda_{ij}) \sum_{s \in \mathcal{S}}\sum_{k \in \mathcal{K}} x_{ijsk} \enspace = \enspace 0 && \forall (i,j) \in \mathcal{A} && \label{eq:nlcon_Uop}\\
& (-U^h_{ij} - \lambda_{ij}) z^h_{ij} \enspace = \enspace 0 && \forall (i,j) \in \mathcal{A} && \label{eq:nlcon_Uh}
\end{flalign}

Note that constraints~\eqref{eq:nlcon_Uop} and~\eqref{eq:nlcon_Uh} are non-linear. To address these, the big M technique is used and binary variables are introduced:
$y^\mathrm{I}_{ij}$ and $y^\mathrm{II}_{ij}$ for constraint~\eqref{eq:nlcon_Uop}; $y^{\mathrm{I}h}_{ij}$ and $y^{\mathrm{II}h}_{ij}$ for constraint~\eqref{eq:nlcon_Uh}.

The only remaining non-linear expression is the first term of the operator's objective function~\eqref{eq:obj_up}. To remedy it, the strong duality theorem can be applied to the lower level problem~\eqref{eq:obj_low}-\eqref{eq:z_con}, as in the work of~\citet{tawfik2019bilevel}. At optimality, we have:
\begin{flalign}
& -\sum_{(i,j) \in \mathcal{A}} \left( \sum_{s \in \mathcal{S}}\sum_{k \in \mathcal{K}} U^O_{ij}x_{ijsk} + \sum_{h \in \mathcal{H}} U^h_{ij}z^h_{ij} \right) \enspace = \enspace \sum_{(i,j) \in \mathcal{A}} D_{ij} \lambda_{ij}  && \label{eq:duality}
\end{flalign}

In addition, the following form is considered for the utility function of the transport operator: 
\begin{flalign}
& U^O_{ij} \enspace = \enspace \bar{U}^O_{ij} + \beta_c p_{ij} + \beta_f f_{ij} \enspace = \enspace \bar{U}^O_{ij} + \beta_c p_{ij} + \beta_f \sum_{s \in \mathcal{S}}\sum_{k \in \mathcal{K}} \phi_{ijs} f_{sk} && \label{eq:u_form}
\end{flalign}
where $\bar{U}_{ij}$ is the part of utility depending on attributes exogenous to the model, $\beta_c$ and $\beta_f$ are the coefficients respectively weighting the importance of price and frequency in the utility function, and $\phi_{ijs}$ is a dummy equal to one if both terminals $i$ and $j$ are contained in service~$s$ and zero otherwise. Then, using equations~\eqref{eq:duality} and~\eqref{eq:u_form}, the first term in~\eqref{eq:obj_up} can be expressed as:
\begin{flalign}
& \sum_{s \in \mathcal{S}}\sum_{k \in \mathcal{K}} p_{ij}x_{ijsk} \enspace = \enspace -\frac{1}{\beta_c} \left( D_{ij} \lambda_{ij}  + \sum_{h \in \mathcal{H}} U^h_{ij}z^h_{ij} + \sum_{s \in \mathcal{S}}\sum_{k \in \mathcal{K}} \bar{U}^O_{ij}x_{ijsk} + \beta_f \sum_{s \in \mathcal{S}}\sum_{k \in \mathcal{K}} \phi_{ijs} f_{sk}x_{ijsk} \right) && \label{eq:ref_px}
\end{flalign}

Because we now have the term $f_{sk}x_{ijsk}$, it may look like a non-linear expression has just been replaced by another. This new term is nevertheless more convenient as the order of magnitude of the frequency is more limited than that of the price. The frequency can then be represented in base~2 conveniently: $f_{sk} = \sum_{b = 0}^{B_{sk}-1} 2^b f_{skb}$, where $f_{skb}$ are binary variables and $B_{sk} = \lceil \log_2(W_{sk}V_k+1) \rceil$. The product term in~\eqref{eq:ref_px} can ultimately be linearized using the well-known technique for the product of binary and continuous variables. The variable representing the product is referred to as $a_{ijskb}$.

The final MILP is then formulated as follows:
\begin{equation}
\begin{split}
\mathbf{(MILP)} \max_{v,f,p,x,z} \quad \sum_{(i,j) \in \mathcal{A}}-\frac{1}{\beta_c} \left( D_{ij} \lambda_{ij}  + \sum_{h \in \mathcal{H}} U^h_{ij}z^h_{ij} + \sum_{s \in \mathcal{S}}\sum_{k \in \mathcal{K}} \bar{U}^O_{ij}x_{ijsk} + \beta_f \sum_{s \in \mathcal{S}}\sum_{k \in \mathcal{K}} \phi_{ijs} \sum_{b = 0}^{B_{sk}-1} 2^b a_{ijskb} \right) \\ - \sum_{s \in \mathcal{S}}\sum_{k \in \mathcal{K}} c^{\mathrm{FIX}}_{sk}f_{sk} - \sum_{(i,j) \in \mathcal{A}}\sum_{s \in \mathcal{S}}\sum_{k \in \mathcal{K}} c^{\mathrm{VAR}}_{ijsk}x_{ijsk} \label{eq:obj_milp}
\end{split}
\end{equation}
\begin{flalign*}
& \mathrm{s.t.} \quad \mathrm{constraints~\eqref{eq:fleet_con}~-~\eqref{eq:vf_con}~\&~\eqref{eq:demand_con}~-~\eqref{eq:z_con}} &&
\end{flalign*}
\begin{flalign}
& f_{sk} \enspace = \enspace \sum_{b = 0}^{B_{sk}-1} 2^b f_{skb} && \forall s \in \mathcal{S}, \enspace \forall k \in \mathcal{K} && \label{eq:fbin_con2}\\
& a_{ijskb} \enspace \le \enspace D_{ij} f_{skb} && \forall (i,j) \in \mathcal{A}, \enspace \forall s \in \mathcal{S}, \enspace \forall k \in \mathcal{K}, \enspace \forall b \in \mathcal{B} && \label{eq:aD_con2}\\
& a_{ijskb} \enspace \le \enspace x_{ijsk} && \forall (i,j) \in \mathcal{A}, \enspace \forall s \in \mathcal{S}, \enspace \forall k \in \mathcal{K}, \enspace \forall b \in \mathcal{B} && \label{eq:ax_con2}\\
& a_{ijskb} \enspace \ge \enspace x_{ijsk} - D_{ij}(1-f_{skb}) && \forall (i,j) \in \mathcal{A}, \enspace \forall s \in \mathcal{S}, \enspace \forall k \in \mathcal{K}, \enspace \forall b \in \mathcal{B} && \label{eq:axD_con2}\\
& \lambda_{ij} \enspace \le \enspace - (\bar{U}^O_{ij} + \beta_c p_{ij} + \beta_f \sum_{s \in \mathcal{S}}\sum_{k \in \mathcal{K}} \phi_{ijs} f_{sk}) && \forall (i,j) \in \mathcal{A} && \label{eq:lUop_con2}\\
& -(\bar{U}^O_{ij} + \beta_c p_{ij} + \beta_f \sum_{s \in \mathcal{S}}\sum_{k \in \mathcal{K}} \phi_{ijs} f_{sk}) - \lambda_{ij} \enspace \le \enspace M^\mathrm{I}_{ij} y^\mathrm{I}_{ij} && \forall (i,j) \in \mathcal{A} && \label{eq:yI_con2}\\
& \sum_{s \in \mathcal{S}}\sum_{k \in \mathcal{K}} x_{ijsk} \enspace \le \enspace M^\mathrm{II}_{ij} y^\mathrm{II}_{ij} && \forall (i,j) \in \mathcal{A} && \label{eq:yII_con2}\\
& y^\mathrm{I}_{ij} + y^\mathrm{II}_{ij} \enspace \le \enspace 1 && \forall (i,j) \in \mathcal{A} && \label{eq:yy_con2}\\
& \lambda_{ij} \enspace \le \enspace - U^h_{ij} && \forall (i,j) \in \mathcal{A}, \enspace \forall h \in \mathcal{H} && \label{eq:lUh_con2}\\
& -U^h_{ij} - \lambda_{ij} \enspace \le \enspace M^{\mathrm{I}h}_{ij} y^{\mathrm{I}h}_{ij} && \forall (i,j) \in \mathcal{A}, \enspace \forall h \in \mathcal{H} && \label{eq:yIh_con2}\\
& z^h_{ij} \enspace \le \enspace M^{\mathrm{II}h}_{ij} y^{\mathrm{II}h}_{ij} && \forall (i,j) \in \mathcal{A}, \enspace \forall h \in \mathcal{H} && \label{eq:yIIh_con2}\\
& y^{\mathrm{I}h}_{ij} + y^{\mathrm{II}h}_{ij} \enspace \le \enspace 1 && \forall (i,j) \in \mathcal{A}, \enspace \forall h \in \mathcal{H} && \label{eq:yyh_con2}\\
& f_{skb} \in \{0,1\} && \forall s \in \mathcal{S}, \enspace \forall k \in \mathcal{K}, \enspace \forall b \in \mathcal{B} && \label{eq:fb_con2}\\
& a_{ijskb} \in \mathbb{N} && \forall (i,j) \in \mathcal{A}, \enspace \forall s \in \mathcal{S}, \enspace \forall k \in \mathcal{K}, \enspace \forall b \in \mathcal{B} && \label{eq:a_con2} \\
& y^\mathrm{I}_{ij}, y^\mathrm{II}_{ij}, y^{\mathrm{I}h}_{ij}, y^{\mathrm{II}h}_{ij} \in \{0,1\} && \forall (i,j) \in \mathcal{A} && \label{eq:y_con2} \\
& \lambda_{ij} \ge 0 && \forall (i,j) \in \mathcal{A} && \label{eq:lam_con2}
\end{flalign}

\subsection{Stochastic formulation}\label{SAA}

In equation~\eqref{eq:u_form}, we set the generic form of the utility function $U^O_{ij}$. In particular, it was assumed that it was fully deterministic, but in reality this is not the case. Firstly, the utility traditionally contains a random term $\epsilon$ (also called \say{error term}), representing the unobserved attributes playing a role in the mode choice. Secondly, one or several $\beta$ coefficients can be assumed as randomly distributed, to account for heterogeneous preferences. Note that these remarks also hold for $U^h_{ij}$. 

With these considerations, the CD-SNDP model becomes a stochastic optimization problem. We then come up with a SAA formulation to solve it. A sample $\mathcal{R}$ is made of $R$ independent realizations. To generate a realization $r$, a draw is performed in the distribution of each random parameter and the corresponding utility functions are computed. For each realization, different utility functions $U^O_{ijr}$ and $U^h_{ijr}$ are obtained: this impacts the mode choice such that the shippers' decision variables $x_{ijskr}$ and $z^h_{ijr}$ become dependent on the sampling. This is also true for the dual variables $\lambda_{ijr}$. Nevertheless, the decision variables of the transport operator are fixed once and for all independently of the sampling. Finally, the objective function is modified such that the average profits over all samples is maximized.


\subsection{Predetermination heuristic}\label{heuristic}

To speed up the solving time of the stochastic formulations, we propose a \say{predetermination heuristic}. As its name suggests, it consists in determining the operator's utility based on given price and frequency values before the optimization. To compute the operator's utility, discrete sets of predefined prices $\mathcal{P}$ and frequencies $\mathcal{F}$ are considered. It is also assumed that the sampling of the shippers' population is already performed so that the utilities of competing alternatives $U^h_{ijr}$ can be computed. Along with the predefined prices $p$ and frequencies $\psi$, these allow to pre-compute several values for each OD pair: the demand faced by the operator $d^{\psi p}_{ij}$, the resulting profit $\pi^{\psi p}_{ij}$, and the price generating the most profit for a given frequency $P^\psi_{ij}$. Algorithm~\ref{alg:price_det} shows the steps to obtain these values.

\begin{algorithm}[H]
\For{$(i,j) \in \mathcal{A}$, $r \in \mathcal{R}$}{
 Determine $U'_{ijr} = \max_{h \in \mathcal{H}} U^h_{ijr}$, and $h'_{ijr} = \argmax_{h \in \mathcal{H}} U^h_{ijr}$\;
}
\For{$(i,j) \in \mathcal{A}$, $\psi \in \mathcal{F}$}{
    \For{$p \in \mathcal{P}$}{
        $d^{\psi p}_{ij} = 0$\;
        \For{$r \in \mathcal{R}$}{
            Compute $U^{O}_{ijr}(\psi,p)$ according to \eqref{eq:u_form}\;
            \If{$U^{O}_{ijr}(\psi,p) \geq U'_{ijr}$}{
                $d^{\psi p}_{ij} = d^{\psi p}_{ij}+\frac{D_{ij}}{|\mathcal{R}|}$\;
            }
        }
        $\pi^{\psi p}_{ij} = p d^{\psi p}_{ij} - \psi \hat{c}^{\mathrm{FIX}}_{ij} - d^{\psi p}_{ij} \hat{c}^{\mathrm{VAR}}_{ij}$\; 
    }
    $P^{\psi}_{ij} = \argmax_{p \in \mathcal{P}} \pi^{\psi p}_{ij}$\;
}
\caption{Price determination method}\label{alg:price_det}
\end{algorithm}

To compute the profit, the fixed and variable costs are needed per OD pair $ij$. However, the available cost parameters are expressed per service $s$ and vehicle type $k$. The costs per OD pair then need to be estimated using the following formulas:
\begin{flalign}
& \hat{c}^{\mathrm{FIX}}_{ij} \enspace = \enspace \sum_{s' \in \mathcal{S}'} \phi_{ijs'} \left( \min_{k} \left( c^{\mathrm{FIX}}_{s'k}/2 \right) \right) && \label{eq:fcosts_est}\\
& \hat{c}^{\mathrm{VAR}}_{ij} \enspace = \enspace \sum_{s' \in \mathcal{S}'} \phi_{ijs'} \left( \min_{k} c^{\mathrm{VAR}}_{ijs'k} \right) && \label{eq:vcosts_est}
\end{flalign}

Here, we consider only services with 2 legs: $\mathcal{S}' = \{s \in \mathcal{S} : |\mathcal{L}_{s}|=2\}$ (i.e. direct services between two terminals $i$ and $j$). For the fixed cost estimation $\hat{c}^{\mathrm{FIX}}_{ij}$, we select the fixed cost of the corresponding service for the cheapest vehicle type and divide it by two (to get the cost for only one service leg). Since the variable costs parameters are already expressed per OD pair, we just select the variable cost of the corresponding service for the cheapest vehicle type as our estimation $\hat{c}^{\mathrm{VAR}}_{ij}$.

Once Algorithm~\ref{alg:price_det} has been used to compute demand values $d^{\psi p}_{ij}$ and price values $P^{\psi}_{ij}$, they can then be used as parameters to solve an auxiliary optimization problem (AP). This problem consists in determining, for a given sample $\mathcal{R}$, the optimal frequencies for fixed prices $\tilde{p}_{ij}$:
\begin{flalign}
&  \mathbf{(AP)} \max_{v,f,g,x,z} \quad \frac{1}{|\mathcal{R}|} \sum_{r \in \mathcal{R}} \left( \sum_{(i,j) \in \mathcal{A}}\sum_{s \in \mathcal{S}}\sum_{k \in \mathcal{K}} \tilde{p}_{ij}x_{ijskr} - \sum_{s \in \mathcal{S}}\sum_{k \in \mathcal{K}} c^{\mathrm{FIX}}_{sk}f_{sk} - \sum_{(i,j) \in \mathcal{A}}\sum_{s \in \mathcal{S}}\sum_{k \in \mathcal{K}} c^{\mathrm{VAR}}_{ijsk}x_{ijskr} \right) && \label{eq:obj_aux}
\end{flalign}
\begin{flalign}
& \mathrm{s.t.} \quad \sum_{s \in \mathcal{S}} v_{sk} \enspace \leq \enspace V_k && \forall k \in \mathcal{K} && \label{eq:fleet_conA} \\
& f_{sk} \enspace \le \enspace W_{sk}v_{sk} && \forall s \in \mathcal{S}, \enspace \forall k \in \mathcal{K} && \label{eq:cycle_conA} \\
& \sum_{r \in \mathcal{R}} \sum_{(i,j) \in \mathcal{A}} \delta_{ijl_{s}}\frac{x_{ijskr}}{|\mathcal{R}|} \enspace \le \enspace Q_k f_{sk} && \forall l_s \in \mathcal{L}_s, \enspace \forall s \in \mathcal{S}, \enspace \forall k \in \mathcal{K} && \label{eq:cap_conA} \\
& x_{ijskr} \enspace \le \enspace \sum_{l_s \in \mathcal{L}_s} \delta_{ijl_{s}}D_{ij} && \forall (i,j) \in \mathcal{A}, \enspace \forall s \in \mathcal{S}, \enspace \forall k \in \mathcal{K}, \enspace \forall r \in \mathcal{R} && \label{eq:serv_conA}\\
& \sum_{s \in \mathcal{S}}\sum_{k \in \mathcal{K}} x_{ijskr} + z_{ijr} \enspace = \enspace D_{ij} && \forall (i,j) \in \mathcal{A}, \enspace \forall r \in \mathcal{R} && \label{eq:demand_conA}\\
& \sum_{\psi \in \mathcal{F}} g_{ij\psi} \enspace \leq \enspace 1 && \forall (i,j) \in \mathcal{A} && \label{eq:binsum_conA} \\
& \sum_{s \in \mathcal{S}}\sum_{k \in \mathcal{K}} \phi_{ijs} f_{sk} \enspace = \enspace  \sum_{\psi \in \mathcal{F}} \psi g_{ij\psi} && \forall (i,j) \in \mathcal{A} && \label{eq:binf_conA} \\
& \sum_{r \in \mathcal{R}}\sum_{s \in \mathcal{S}}\sum_{k \in \mathcal{K}} \frac{x_{ijskr}}{|\mathcal{R}|} \enspace \leq \enspace \sum_{\psi \in \mathcal{F}} g_{ij\psi} d^{\psi \tilde{p}}_{ij} && \forall (i,j) \in \mathcal{A} &&  \label{eq:dembin_conA} \\
& v_{sk} \in \mathbb{N} && \forall s \in \mathcal{S}, \enspace \forall k \in \mathcal{K} && \label{eq:v_conA}\\
& f_{sk} \in \mathbb{N} && \forall s \in \mathcal{S}, \enspace \forall k \in \mathcal{K} && \label{eq:f_conA}\\
& g_{ij\psi} \in \{0,1\} && \forall (i,j) \in \mathcal{A}, \enspace \forall s \in \mathcal{S}, \enspace \forall \psi \in \mathcal{F} && \label{eq:g_conA}\\
& x_{ijskr} \ge 0 && \forall (i,j) \in \mathcal{A}, \enspace \forall s \in \mathcal{S}, \enspace \forall k \in \mathcal{K}, \enspace \forall r \in \mathcal{R} && \label{eq:x_conA}\\
& z_{ijr} \ge 0 && \forall (i,j) \in \mathcal{A}, \enspace \forall r \in \mathcal{R} && \label{eq:z_conA}
\end{flalign}

This auxiliary problem contains additional elements that deserve some discussion. First, the objective~\eqref{eq:obj_aux} is now formulated as a SAA function and the decision variables $x$ and $z$ are now dependent on $r$. Constraints~\eqref{eq:cap_conA} to~\eqref{eq:demand_conA} are modified accordingly. A new binary variable $g_{ij\psi}$ is introduced: it is equal to one if the predefined frequency $\psi$ is chosen for OD pair $(i,j)$, and zero otherwise. Constraint~\eqref{eq:binsum_conA} ensures that at most one frequency $\psi$ is chosen per OD pair. The value of $\psi$ is then linked to the decision variable of services frequency $f$ through constraint~\eqref{eq:binf_conA}. Finally, constraint~\eqref{eq:dembin_conA} aggregates the decision variables $x_{ijskr}$ of cargo assigned to the operator over the whole sample and bounds it with the precomputed demand $d^{\psi p}_{ij}$ defined in Algorithm~\ref{alg:heur}. This last constraint allows to keep the utility functions out of the optimization problem. As a result, the variable $z_{ijr}$ is now independent of the competing modes $h$. Once the optimization is performed, the corresponding value of $z^*_{ijr}$ can be assigned to the competing mode $h'_{ijr}$ with the maximum utility as computed in Algorithm~\ref{alg:price_det}.

Getting rid of the utilities and the pricing decision in the optimization allows to considerably decrease the solving time. Indeed, the variables $p_{ij}$, $f_{skb}$, $a_{ijskb}$, $\lambda_{ij}$ and $y_{ij}$'s are not used anymore, and only the variables $g_{ij\psi}$ are added. The number of constraints is also drastically reduced. The idea of the heuristic is to exploit this advantage to solve the auxiliary problem iteratively, as described in Algorithm~\ref{alg:heur}.

\begin{algorithm}[H]
Use Algorithm~\ref{alg:price_det} to pre-compute $d^{\psi p}_{ij}$ and $P^{\psi}_{ij}$\;
Define the set of visited solutions $\Omega = \emptyset$\;
Set each $\tilde{p}_{ij}$ with arbitrary prices contained in $\mathcal{P}$\;
Set each $\tilde{\psi}_{ij}$ with arbitrary frequencies contained in $\mathcal{F}$\;
\While{$(\tilde{p}_{ij},\tilde{\psi}_{ij}) \notin \Omega$}{
    Solve (AP) to get $g^*_{ij\psi}$, i.e. the chosen frequency $\psi$ corresponding to predefined prices $\tilde{p}_{ij}$\;
    Update each $\tilde{\psi}_{ij}$ with the frequency associated to $g^*_{ij\psi}$\;
    Update each $\tilde{p}_{ij}$ with the value $P^{\tilde{\psi}_{ij}}_{ij}$\;
    Add $(\tilde{p}_{ij},\tilde{\psi}_{ij})$ to the set $\Omega$\;
}
Return the solution $(\tilde{p}_{ij},\tilde{\psi}_{ij})$\;
\caption{Predetermination heuristic}\label{alg:heur}
\end{algorithm}


The performance of the heuristic is highly dependent on the size of sets $\mathcal{P}$ and $\mathcal{F}$. The more values they contain, the better is the approximation at the cost of additional computational resources. These sets should then ensure a good coverage of the search space in order for the heuristic to return satisfying solutions.


\section{Case Study}\label{case}

The proposed CD-SNDP is applied to container transport on a small network of 3 nodes: Rotterdam (RTM), Duisburg (DUI) and Bonn (BON). We consider an inland vessel operator competing with two other modes (Road and Rail) and another IWT carrier.

\subsection{Description}\label{description}

The operator's fleet is composed of 2 vessel types : 24 vessels of type M8 and 12 vessels of type M11 with a maximal capacity of 180 TEUs and 300 TEUs, respectively. Each vessel type has a maximal operation time, $T^{max}$, of 120 hours per week. The transport demand inputs are based on the NOVIMOVE project~\citep{D_2_2}, whereas the costs for IWT and the two competing modes as well as the sailing time $t^{sail}$ and the time spent in ports $t^{port}$ are estimated using the model of~\citet{shobayo2021conceptual}. Based on these inputs, the computation of the maximum number of cycles is straightforward: ${W_{sk}=T^{max}/(t^{sail}_{sk}+t^{port}_{sk})}$.

Regarding shippers, the utility functions are formulated as follows:
\begin{flalign}
    & U^O_{ij} & = & \alpha^{\mathrm{IWT}} + \beta_{a}^\mathrm{Inter} a^{\mathrm{IWT}}_{ij} + \beta_{q}^\mathrm{IWT} q_{ij} + \beta_{c}^\mathrm{Inter} (\pmb{p_{ij}}+\mathrm{VoT}t^{\mathrm{IWT}}_{ij}) + \beta_{f}^\mathrm{Inter} \pmb{f_{ij}} + \epsilon^O_{ij} && \label{eq:UO} \\
   & U_{ij}^{h=\mathrm{IWT}} & = & \alpha^{\mathrm{IWT}} + \beta_{a}^\mathrm{Inter} a^{\mathrm{IWT}}_{ij} + \beta_{q}^\mathrm{IWT} q_{ij} + \beta_{c}^\mathrm{Inter} (p^{\mathrm{IWT}}_{ij}+\mathrm{VoT}t^{\mathrm{IWT}}_{ij}) + \beta_{f}^\mathrm{Inter} f^{\mathrm{IWT}}_{ij} + \epsilon^{\mathrm{IWT}}_{ij} && \label{eq:UIWT} \\
   & U_{ij}^{h=\mathrm{Rail}} & = & \alpha^{\mathrm{Rail}} + \beta_{a}^\mathrm{Inter} a_{ij}^{\mathrm{Rail}} + \beta_{c}^\mathrm{Inter} (p^{\mathrm{Rail}}_{ij}+\mathrm{VoT}t^{\mathrm{Rail}}_{ij}) + \beta_{f}^\mathrm{Inter} f_{ij}^{\mathrm{Rail}} + \epsilon^{\mathrm{Rail}}_{ij} && \label{eq:URail} \\
   & U_{ij}^{h=\mathrm{Road}} & = & \alpha^{\mathrm{Road}} + \beta_{a}^\mathrm{Road} a_{ij}^{\mathrm{Road}} + \beta_{c}^\mathrm{Road} (p^{\mathrm{Road}}_{ij}+\mathrm{VoT}t^{\mathrm{Road}}_{ij}) + \epsilon^{\mathrm{Road}}_{ij} && \label{eq:URoad}
\end{flalign}
where, for each mode, $\alpha$ is the alternative specific constant, $a$ is an accessibility metric, $t$ is the total travel time in hours, $p$ is the cost for shippers in thousand of euros per TEU, and $f$ is the weekly frequency for intermodal transports (i.e. IWT and Rail). Moreover, $q_{ij}$ is a dummy equal to one if a seaport is located at $i$ or $j$ and $\mathrm{VoT}$ is the Value of Time expressed in 1000\euro{}/TEU/hour. Each attribute is weighted by a coefficient $\beta$ and each mode has a random error term $\epsilon$. Although they have similarities, it is assumed that the vessel operator and the IWT carrier alternatives are not correlated. The same assumption holds between all alternatives. Therefore, in the remainder of this work, all the error terms $\epsilon_{ij}$ are considered independent and identically distributed (iid).

Within the CD-SNDP context, all the terms contained in the utilities of the competing modes (IWT, Rail and Road) are exogenous to the optimization model and are thus treated as parameters. Regarding the utility of the operator, only the terms in bold ($p_{ij}$ and $f_{ij}$) are endogenous while the other terms are also parameters. $p_{ij}$ is the decision variable on pricing and $f_{ij}$ corresponds to the term $\sum_{s \in \mathcal{S}}\sum_{k \in \mathcal{K}} \phi_{ijs} f_{sk}$, as introduced in equation~\eqref{eq:u_form}.

The model's coefficients were estimated with aggregate data using a Weighted Logit methodology. It is named \say{weighted} because, during the estimation, the log-likelihood function is weighted by the yearly cargo flows on each OD pair~\citep{Rich2009}. It thus gives more importance to the OD pairs with high volumes. For more details, the reader is referred to~\citet{nicolet2022logit}. One noteworthy characteristic of the data on which the coefficients were estimated is that the frequency for IWT does not exceed 35 services per week. Therefore, the following constraint is added to our CD-SNDP problem to guarantee consistency between the results and the mode choice model:
\begin{flalign}
& f_{sk} \enspace \leq \enspace 35 \qquad \forall s \in \mathcal{S}, \enspace \forall k \in \mathcal{K} && \label{eq:bound_freq}
\end{flalign}

We use this case study to compare the results of 3 deterministic and 2 stochastic models. The former consist in the benchmark, the SNDP and the CD-SNDP, which uses only the deterministic part of the utility functions in equations \eqref{eq:UO} to \eqref{eq:URoad} without error terms $\epsilon$. The latter two are stochastic variations of the CD-SNDP:
\begin{itemize}
    \item Multinomial Logit (MNL): with iid error terms $\epsilon$, following an Extreme Value distribution;
    \item Mixed Logit: random $\beta_{c}^\mathrm{Inter}$ following a Lognormal distribution of parameters $\mu_c^{\mathrm{Inter}}$ and $\sigma_c^{\mathrm{Inter}}$ (representing the heterogeneous cost sensitivity of shippers) together with iid error terms $\epsilon$.
\end{itemize}

\subsection{Evaluation through out-of-sample simulation}\label{simul}

In order to assess the solutions returned by these models, we simulate the demand response using an out-of-sample population. Indeed, the profit returned by the optimization is the one expected based on the SAA and the model's assumptions, but it gives no indication on how well the solution will perform with actual shippers. This out-of-sample simulation also allows to compare the different models with each other. The procedure is as follows:
\begin{enumerate}
    \item For each OD pair, generate a population of 1000 shippers (i.e. perform 1000 draws of $\epsilon$ and $\beta_{c}^\mathrm{Inter}$, note that these draws are different than the ones used in the SAA) and divide the demand $D_{ij}$ equally among the shippers;
    \item For each shipper, compute their utilities by plugging the drawn $\epsilon$ and $\beta_{c}^\mathrm{Inter}$, as well as the frequencies and prices returned by the model, into equations \eqref{eq:UO} to \eqref{eq:URoad};
    \item Allocate the shipper's containers to the alternative with the maximal utility;
    \item When all containers have been allocated, compute the resulting modal shares and the actual profit for the inland vessel operator.
\end{enumerate}

\subsection{Coefficients of utility functions}\label{uf_coeff}

For the out-of-sample simulation, we directly make use the coefficients of the Weighted Logit Mixture model estimated in~\citet{nicolet2022logit}. However, these true utility functions of the shippers are not known by the operator. The same coefficients cannot, therefore, be used in the CD-SNDP. To alleviate this issue, we use the Weighted Logit Mixture to generate synthetic choice data, from which utility coefficients can be estimated by the operator. This process ensures that the true utility functions remain hidden from the operator, as they only have access to the choice realizations of shippers.

The available inputs are the OD matrices and the attributes related to IWT, Rail and Road on each OD pair along the European Rhine-Alpine (RA) corridor. To generate a choice instance for a given OD pair using the Weighted Logit Mixture, we first draw the value of $\beta_{c}^\mathrm{Inter}$ and each mode's $\epsilon$ in their respective distributions. Then, they are plugged, along with each mode's attributes, into equations \eqref{eq:UIWT} to \eqref{eq:URoad}. Finally, the mode with the highest utility is selected and we get one synthetic choice instance. This process is then repeated for all OD pairs. To remain consistent with the Weighted Logit methodology, the number of generated choice instances per OD pair is set proportional to its cargo volume. In particular, each OD pair get at least one choice instance and an additional instance is generated per 10'000 TEUs circulating yearly on the OD pair. As a result, we end up with a synthetic dataset composed of 8676 choice instances, from which the MNL and Mixed Logit models can be estimated.

The coefficients of the Weighted Logit Mixture model are presented in Table~\ref{tab:coeffU}, along with the coefficients of the Mixed Logit and MNL estimated using the synthetic dataset (note that $\alpha_{\mathrm{IWT}}$ is normalized to zero). The mean value of $\beta_{c}^\mathrm{Inter}$ is also presented.

\begin{table}[!t]
\caption{Coefficients of the mode choice models}
\label{tab:coeffU}
\setlength{\tabcolsep}{10pt}
\centering
\begin{tabular}{c|c||c|c}
 & \begin{tabular}[c]{@{}c@{}}Actual population\end{tabular} & \multicolumn{2}{c}{Synthetic data} \\ 
\hline
Parameter & \begin{tabular}[c]{@{}c@{}}Weighted Logit\\ Mixture\end{tabular} & \begin{tabular}[c]{@{}c@{}}Mixed\\ Logit\end{tabular} & MNL \\ \hline
\rule{0pt}{2.6ex}
$\alpha^{\mathrm{IWT}}$ & 0 & 0 & 0 \\
$\alpha^{\mathrm{Rail}}$ & 0.713 & 0.816 & 0.338 \\
$\alpha^{\mathrm{Road}}$ & 2.30 & 2.35 & 2.06 \\
$\beta_q^{\mathrm{IWT}} $ & 1.63 & 1.60 & 1.49 \\
$\beta_f^{\mathrm{Inter}}$ & 0.0278 & 0.0262 & 0.0229 \\
$\beta_a^{\mathrm{Road}}$ & 0.0530 & 0.0506 & 0.0469 \\
$\beta_a^{\mathrm{Inter}}$ & 0.157 & 0.173 & 0.141 \\
$\beta_c^{\mathrm{Road}}$ & -8.68 & -8.73 & -4.81 \\
$\beta_c^{\mathrm{Inter}}$ & \cellcolor[HTML]{EFEFEF} & \cellcolor[HTML]{EFEFEF} & -5.76 \\ \hline
\rule{0pt}{2.6ex}
$\mu_c^{\mathrm{Inter}}$ & 2.30 & 2.40 & \cellcolor[HTML]{EFEFEF}  \\
$\sigma_c^{\mathrm{Inter}}$ & 0.690 & 0.618 & \cellcolor[HTML]{EFEFEF} \\
$\bar{\beta}_c^{\mathrm{Inter}}$ & -12.65 & -13.34 & -5.76 \\ \hline

\end{tabular}
\end{table}

\subsection{Deterministic results}\label{results}

In the remainder of this section, we will present and discuss the results of these various models, starting with the deterministic ones.

\subsubsection{CD-SNDP vs. Benchmark and SNDP}\label{res_deter}

The weekly frequencies for both vessel types and the charged prices are reported in Table~\ref{tab:resdeter}. In order to better understand the pricing decision, the table also displays the prices of the competing alternatives. For the benchmark and SNDP, the optimal prices are set at the same level as the cheapest competing alternative (in our case, IWT). This is because of the assumption that shippers are purely cost-minimizers and the deterministic nature of the models: if the vessel operator charges just 0.001~\euro{} less than the cheapest alternative, then the models will consider that all shippers will choose the services of the vessel operator instead of the competition. In the CD-SNDP, shippers are assumed to consider other attributes beside cost to perform their mode choice: the optimal prices then differ from the cheapest alternative.

\begin{table}[!t]
\centering
\caption{Solutions of deterministic models with prices of the competing alternatives}
\label{tab:resdeter}
\resizebox{0.9\textwidth}{!}{%
\begin{tabular}{llccc|ccc}
\multicolumn{1}{c}{} & \multicolumn{1}{c}{} & Benchmark & SNDP & CD-SNDP & IWT & Road & Rail \\ \hline
\multirow{6}{*}{Prices {[}\euro{}{]}} & RTM-DUI & 68 & 68 & 148 & 68 & 252 & 203 \\
 & DUI-RTM & 69 & 69 & 160 & 69 & 251 & 203 \\
 & RTM-BON & 76 & 76 & 80 & 76 & 317 & 214 \\
 & BON-RTM & 74 & 74 & 51 & 74 & 315 & 214 \\
 & DUI-BON & 46 & 46 & - & 46 & 136 & 152 \\
 & BON-DUI & 46 & 46 & - & 46 & 136 & 152 \\ \hline
 \multirow{4}{*}{\begin{tabular}[c]{@{}l@{}}Weekly\\ frequencies\\ {[}M8 vessels \\ (M11 vessels){]}\end{tabular}} & RTM-DUI & 16 (12) & 0 (13) & 0 (18) & \cellcolor[HTML]{EFEFEF} & \cellcolor[HTML]{EFEFEF} & \cellcolor[HTML]{EFEFEF} \\
 & RTM-BON & 0 (5) & 0 (5) & 0 (0) & \cellcolor[HTML]{EFEFEF} & \cellcolor[HTML]{EFEFEF} & \cellcolor[HTML]{EFEFEF} \\
 & DUI-BON & 32 (3) & 20 (0) & 0 (0) & \cellcolor[HTML]{EFEFEF} & \cellcolor[HTML]{EFEFEF} & \cellcolor[HTML]{EFEFEF} \\
 & RTM-DUI-BON & - & 19 (0) & 24 (0) & \cellcolor[HTML]{EFEFEF} & \cellcolor[HTML]{EFEFEF} & \cellcolor[HTML]{EFEFEF} \\ \hline
\end{tabular}%
}
\end{table}

\begin{figure}[!b]
\centering
\includegraphics[width=0.8\textwidth]{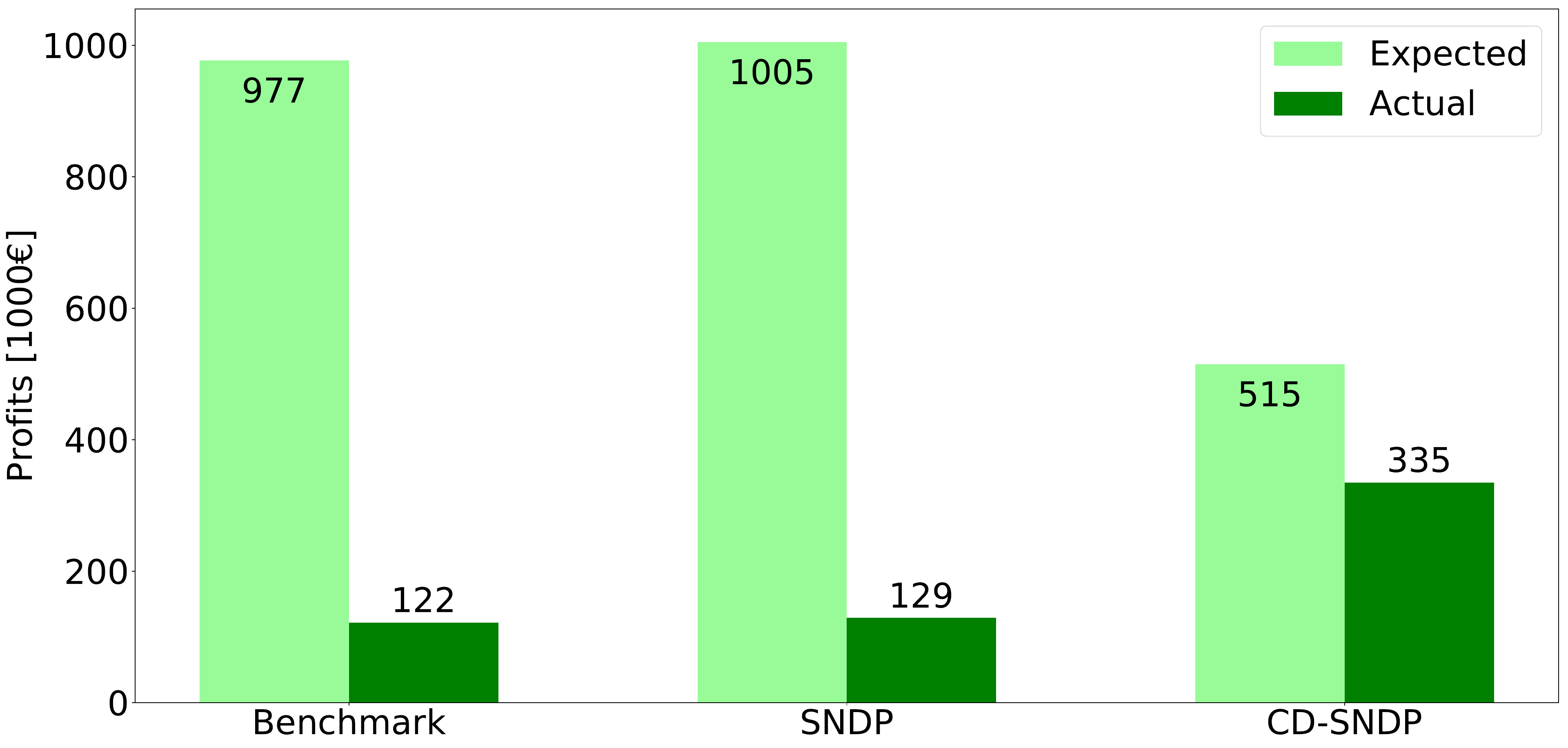}
\caption{Comparison of profits obtained with the deterministic models.\label{fig:DeterProfits}}
\end{figure}
\begin{figure}[!t]
\centering
\includegraphics[width=0.8\textwidth]{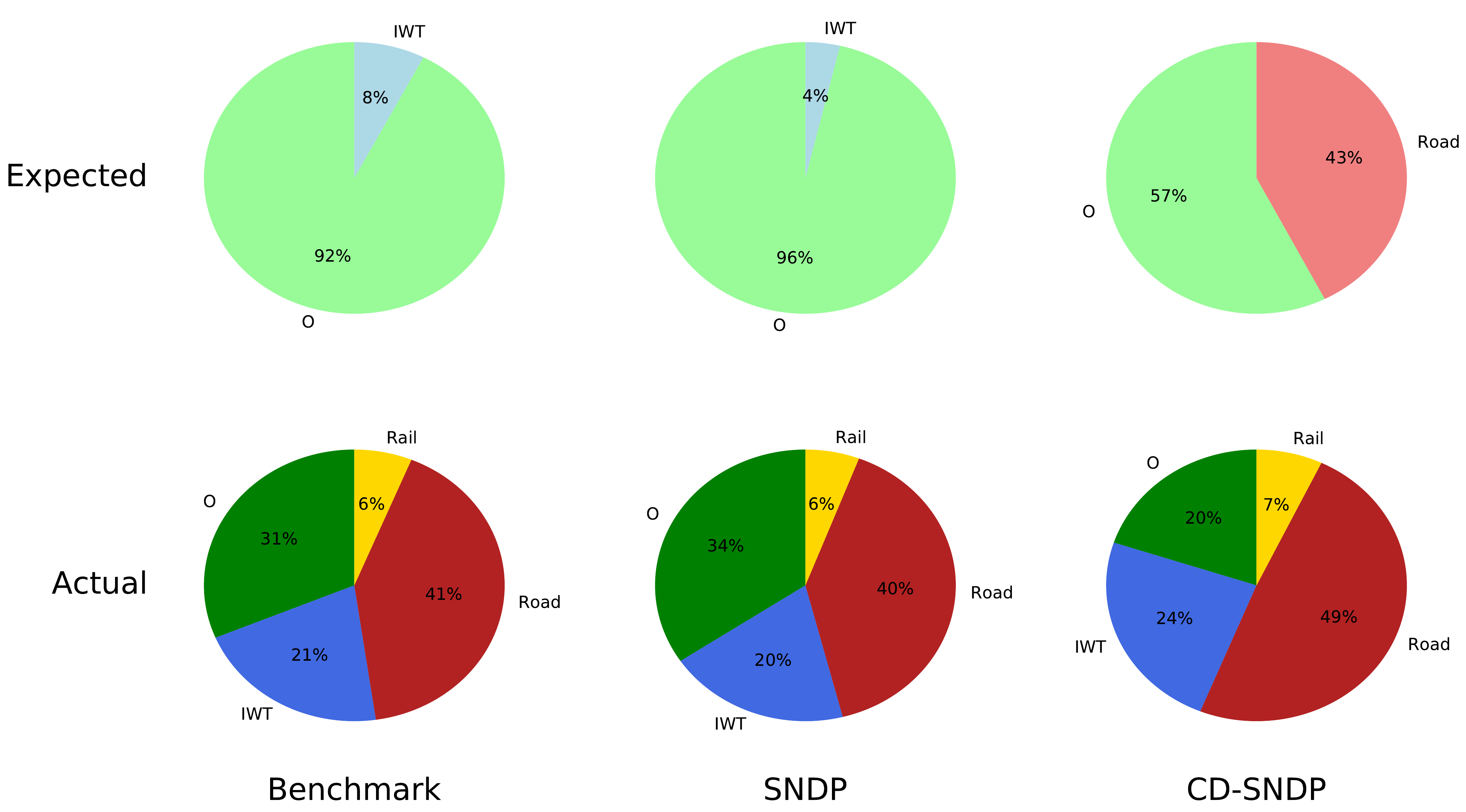}
\caption{Comparison of modal shares returned by the deterministic models and simulated ones.\label{fig:DeterShares}}
\end{figure}

Regarding the optimal frequencies, allowing to visit more than 2 terminals per service provides additional flexibility to the SNDP compared to the Benchmark. The SNDP takes advantage of this consolidation opportunity which results in higher expected profits. Figure~\ref{fig:DeterProfits} displays the expected profits versus the actual ones returned by the out-of-sample simulation: it shows that the SNDP also returns higher simulated profits. The reason is that the vessel operator is able to attract more demand with this 3-stop service. This is seen in Figure~\ref{fig:DeterShares}, which represents the expected and actual modal shares for each deterministic model.

For the CD-SNDP, the expected profits drop compared to the two other models. This is because the OD pair Duisburg-Bonn is not served anymore by the vessel operator. The distance between these two terminals is indeed relatively short so, as the CD-SNDP takes multiple factors in consideration for the mode choice of shippers, it is evident that Road becomes the preferred option for this OD pair. Although the vessel operator gets smaller market shares than with the Benchmark and SNDP (See Figure~\ref{fig:DeterShares}), the CD-SNDP returns actual profits that are more than 2.5 times higher. Since the choice-driven model also considers frequency in the mode choice, it is able to charge much more on the OD pair Rotterdam-Duisburg due the very high proposed frequency (42 services per week) between these two terminals.

These deterministic results already suggest that significant gains can be achieved with the Choice-Driven SNDP. More efficient services and pricing can be designed, thus resulting in considerably increased profits.

\subsection{Stochastic variants}\label{res_stoch}

In this section, the results of the stochastic versions of the CD-SNDP are described. Two random utility formulations are compared: MNL (with random error terms $\epsilon$) and Mixed Logit (with $\epsilon$ and distributed cost sensitivity $\beta_{c}^\mathrm{Inter}$). First, we present the results obtained with the exact method, and then the ones with the predetermination heuristic.

\subsubsection{Exact method}\label{res_stoch_ex}

The two stochastic variants are solved through SAA with a sample size of $R = 1000$, i.e. a thousand draws are performed in the distributions of $\epsilon$ (and of $\beta_{c}^\mathrm{Inter}$, for the Mixed Logit). For each variant, we run 20 replications with 20 different random seeds, thus generating 20 different samples. The aggregated statistics of the obtained solutions and computation time are presented in Table~\ref{tab:resexact}. Note that a time limit of 48 hours has been applied, that is why the statistics of the optimality gap are also presented.

\begin{table}[!t]
\centering
\caption{Solutions of stochastic models with exact method (1000 draws)}
\label{tab:resexact}
\resizebox{0.9\textwidth}{!}{%
\begin{tabular}{llccc|ccc}
\multicolumn{1}{c}{} & \multicolumn{1}{c}{} & \multicolumn{3}{c}{MNL} & \multicolumn{3}{c}{Mixed Logit} \\ \hline
\multicolumn{1}{c}{} & \multicolumn{1}{c}{} & Min. & Average & Max. & Min. & Average & Max. \\ \hline
\multirow{4}{*}{\begin{tabular}[c]{@{}l@{}}Weekly\\ frequencies\\ {[}M8 vessels \\ (M11 vessels){]}\end{tabular}} & RTM-DUI & 0 (0) & 0 (0) & 0 (0) & 0 (0) & 0 (0) & 0 (0) \\
 & RTM-BON & 0 (0) & 0 (0) & 0 (0) & 0 (0) & 0 (0) & 0 (0) \\
 & DUI-BON & 0 (0) & 0 (0) & 0 (0) & 0 (0) & 0 (0) & 0 (0) \\
 & RTM-DUI-BON & 24 (0) & 24 (0) & 24 (0) & 24 (0) & 24 (0) & 24 (0) \\ \hline
\multirow{6}{*}{Prices {[}\euro{}{]}} & RTM-DUI & 160 & 197 & 245 & 132 & 161 & 193 \\
 & DUI-RTM & 160 & 192 & 255 & 120 & 147 & 180 \\
 & RTM-BON & 180 & 205 & 253 & 153 & 194 & 228 \\
 & BON-RTM & 160 & 189 & 245 & 140 & 182 & 225 \\
 & DUI-BON & 106 & 163 & 242 & 98 & 140 & 180 \\
 & BON-DUI & 106 & 163 & 242 & 113 & 146 & 184 \\ \hline
\multicolumn{2}{l}{Computation time {[}hours{]}} & 5.0 & 34.5 & 48.0 & 48.0 & 48.0 & 48.0 \\ \hline
\multicolumn{2}{l}{Optimality gap} & 0.0\% & 6.5\% & 28.8\% & 28.0\% & 84.9\% & 111.8\% \\ \hline
\end{tabular}%
}
\end{table}

The frequency decision remains identical for all 20 replications of both stochastic variants. The pricing decision, on the other hand, is very variable from one replication to another. The variation is slightly more pronounced for the MNL case than for the Mixed Logit, but the main takeaway is that the MNL results in higher prices than the Mixed Logit. Also, both variants find higher prices than the deterministic CD-SNDP.

\begin{figure}[!b]
\centering
\includegraphics[width=0.8\textwidth]{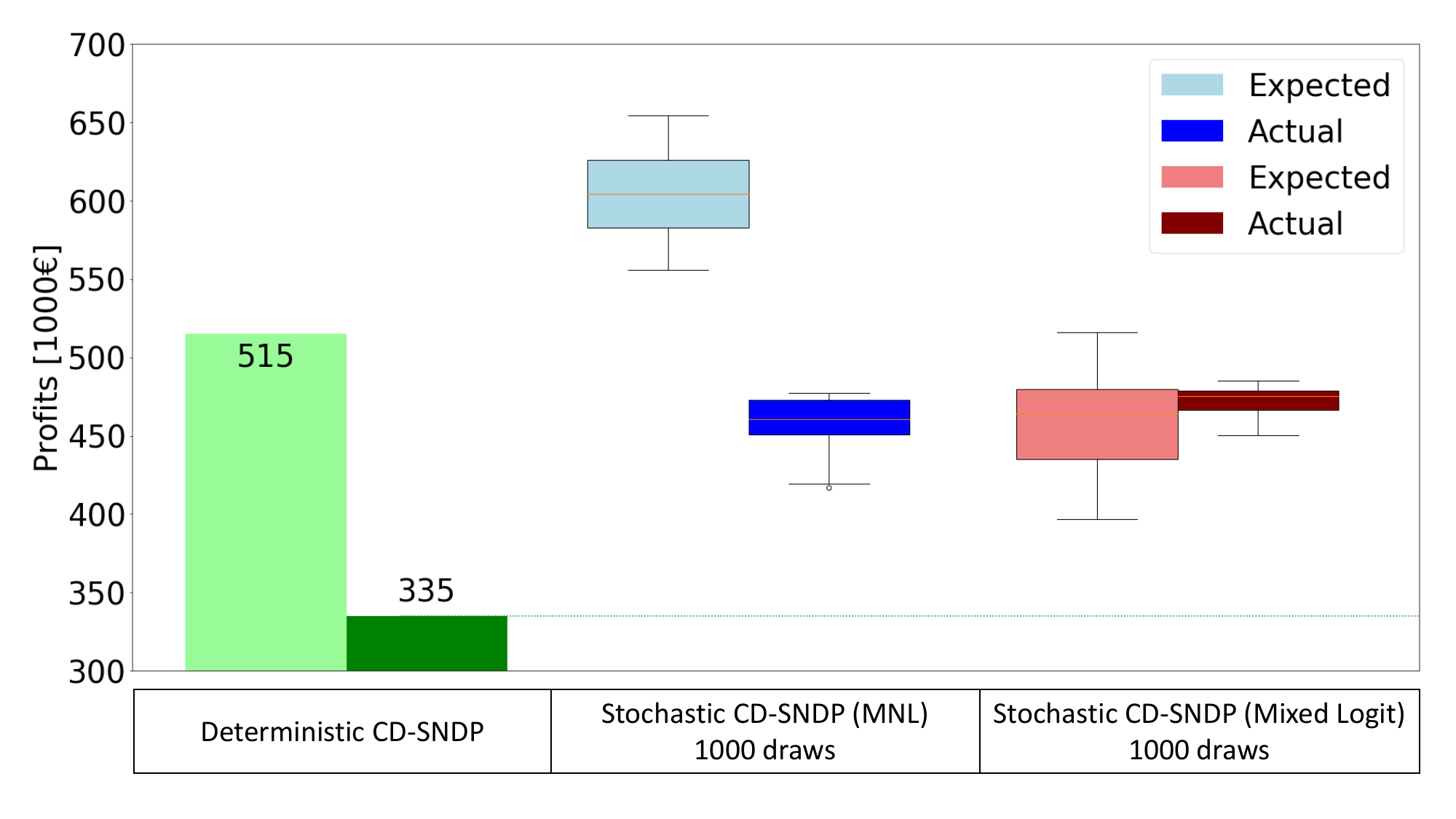}
\caption{Comparison of profits by the stochastic models with exact method and the deterministic CD-SNDP.\label{fig:ExactProfits}}
\end{figure}

\begin{figure}[!t]
\centering
\includegraphics[width=0.8\textwidth]{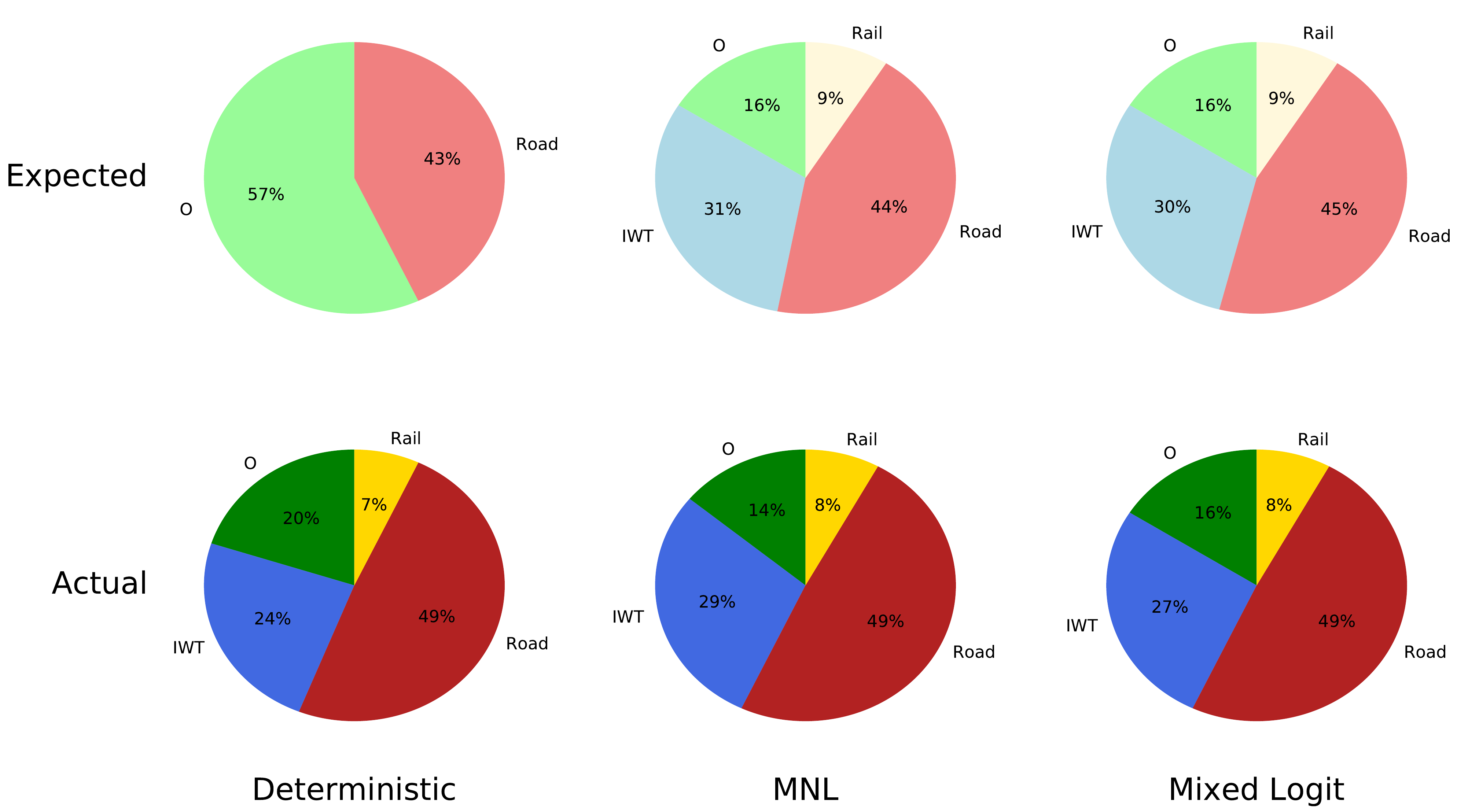}
\caption{Comparison of average modal shares returned by the stochastic models with exact method and simulated ones.\label{fig:ExactShares}}
\end{figure}

The influence on the expected and actual profits is depicted in the boxplots of Figure~\ref{fig:ExactProfits}. The higher prices set by the MNL lead to greater expected profits compared to the Mixed Logit. But this difference is canceled out when comparing the simulated profits as the MNL profits fall at a very slightly lower level than the ones of the Mixed Logit. Nevertheless, the actual profits for both variants are substantially higher than for the deterministic CD-SNDP. In fact, they provide an additional 40\% gain compared to the deterministic case. This is because, due to the more detailed choice models, the modal shares can be estimated much better during the optimization. Indeed, Figure~\ref{fig:ExactShares} shows the average expected modal shares against the actual ones. These shares are really close to each other for both the MNL and the Mixed Logit, whereas the deterministic model highly overestimates the share of the vessel operator.

This also explains why no bigger vessels are being operated in both stochastic variants, compared to the deterministic case. Due to the very high operator share estimated in the latter, smaller vessels only are not sufficient to meet the demand. The stochastic solutions prefer to operate only smaller vessels, as the operating costs are lower and the service capacity is sufficient to accommodate all the incoming demand.

Comparing the MNL with the Mixed Logit, the accuracy of their modal share estimation is nearly equivalent, but the MNL tends to slightly overestimate the operator's share during the optimization process. This is why the expected profits in Figure~\ref{fig:ExactProfits} are significantly higher than the actual ones. On the other hand, the expected profits with the Mixed Logit are in line with the actual ones. In the end, the actual profits for both the MNL and Mixed Logit end up at a similar level. However, the large optimality gaps reported for the Mixed Logit in Table~\ref{tab:resexact} prevent any conclusion at this stage. Indeed, no replication was able to terminate within the 48 hours limit. Even though the addition of stochasticity in the CD-SNDP provide more gains, it is done at the expense of computing time. In order to remedy this, we make use of the predetermination heuristic presented in Section~\ref{heuristic}.

\subsubsection{Predetermination heuristic}\label{res_stoch_heur}

The two stochastic variants are solved using the same samples as for the exact method. We use the following set of predefined prices: $\mathcal{P}=\mathbb{N} \cap [0,500]$, and the set of predefined frequencies: $\mathcal{F}=\mathbb{N} \cap [0,35]$ in accordance with \eqref{eq:bound_freq}. The statistics of the heuristic solutions are reported in Table~\ref{tab:resheur} together with the computation time.

\begin{table}[!b]
\centering
\caption{Solutions of stochastic models with predetermination heuristic (1000 draws)}
\label{tab:resheur}
\resizebox{0.9\textwidth}{!}{%
\begin{tabular}{llccc|ccc}
\multicolumn{1}{c}{} & \multicolumn{1}{c}{} & \multicolumn{3}{c}{MNL} & \multicolumn{3}{c}{Mixed Logit} \\ \hline
\multicolumn{1}{c}{} & \multicolumn{1}{c}{} & Min. & Average & Max. & Min. & Average & Max. \\ \hline
\multirow{4}{*}{\begin{tabular}[c]{@{}l@{}}Weekly\\ frequencies\\ {[}M8 vessels \\ (M11 vessels){]}\end{tabular}} & RTM-DUI & 0 (0) & 0 (8) & 1 (14) & 0 (0) & 0 (4) & 0 (14) \\
 & RTM-BON & 0 (0) & 0 (0) & 0 (0) & 0 (0) & 0 (0) & 0 (0) \\
 & DUI-BON & 0 (0) & 8 (2) & 16 (6) & 0 (0) & 8 (2) & 17 (5) \\
 & RTM-DUI-BON & 16 (0) & 22 (0) & 24 (0) & 18 (0) & 22 (0) & 24 (0) \\ \hline
\multirow{6}{*}{Prices {[}\euro{}{]}} & RTM-DUI & 163 & 198 & 239 & 122 & 157 & 197 \\
 & DUI-RTM & 169 & 201 & 247 & 119 & 150 & 183 \\
 & RTM-BON & 174 & 203 & 249 & 141 & 189 & 225 \\
 & BON-RTM & 166 & 190 & 225 & 134 & 177 & 223 \\
 & DUI-BON & 111 & 170 & 265 & 92 & 159 & 247 \\
 & BON-DUI & 111 & 172 & 265 & 81 & 157 & 247 \\ \hline
\multicolumn{2}{l}{Computation time {[}hours{]}} & 0.08 & 0.09 & 0.11 & 0.10 & 0.11 & 0.12 \\ \hline
\end{tabular}%
}
\end{table}

Compared to the exact method, the predetermination heuristic is remarkably faster. When the computation was in the order of days for the exact method, it is now reduced to a few minutes. Most of these minutes are spent precomputing the demand and price values with Algorithm~\ref{alg:price_det}. With the heuristic, there is also little difference in solving time between the two stochastic variants.

\begin{figure}[!t]
\centering
\includegraphics[width=0.8\textwidth]{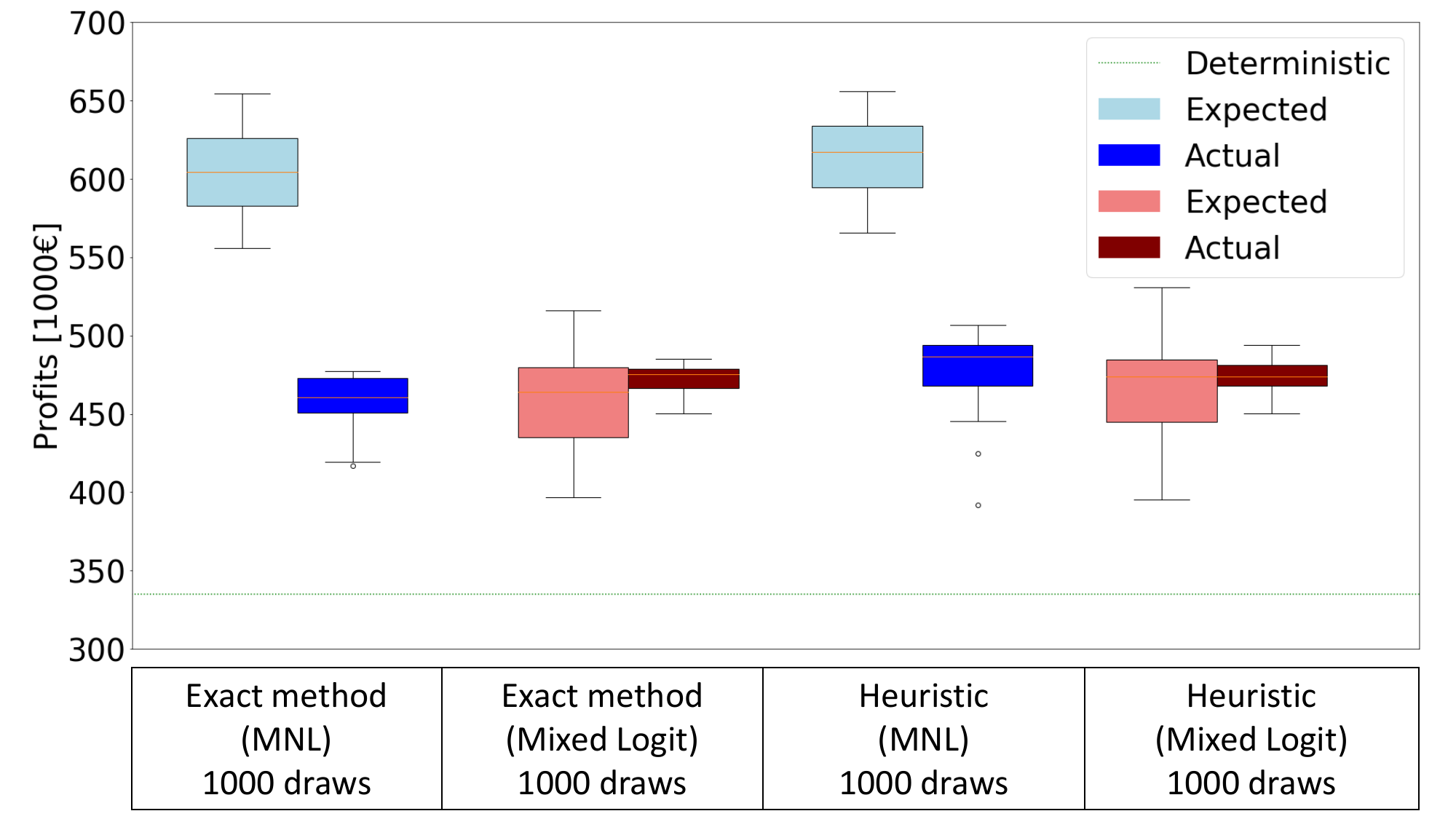}
\caption{Comparison of profits by the stochastic models with exact method and the heuristic.\label{fig:ExactHeuristic}}
\end{figure}

Regarding the quality of the solutions, the prices found with the heuristic are consistent with the ones returned by the exact method. However, there is some variation in the resulting frequencies. With the exact method, the optimal solution for all replications with both choice models was deploying all the small vessels on the 3-stop service and not using the bigger vessels. With the heuristic, small vessels are still concentrated on the 3-stop service but some are also deployed on the service between Duisburg and Bonn. Some bigger vessels are also used on the services between Rotterdam and Duisburg and between Duisburg and Bonn, but not between Rotterdam and Bonn because the precomputed profits are decreasing with the increase of frequency on this OD pair. The solution of the heuristic makes a better use of the operator's fleet: it allows the operator to attract more demand on certain OD pairs due to the increased frequencies, but additional costs are incurred. The comparison between the profits obtained with the exact method and with the predetermination heuristic is shown in Figure~\ref{fig:ExactHeuristic}.

The profit ranges found by the heuristic are similar to the ones with the exact method. Even with increased frequencies (generating additional costs) and similar price ranges, the profits remain stable because the increased frequencies attract more demand, compensating for the incurred costs. We still observe a significant gap between the expected and actual profits in the MNL case, whereas these two values are at a similar level for the Mixed Logit. This is because the MNL used in the CD-SNDP has a much lower cost coefficient $\beta_c^{\mathrm{Inter}}$ in absolute value than in the actual population (see Table~\ref{tab:coeffU}). On the other hand, the Mixed Logit has coefficients that are more in line with the actual population. The cost sensitivity of the shippers is then underestimated by the MNL, which results in prices that are higher than with the Mixed Logit. The CD-SNDP with MNL then expects that high profits will be realized, whereas in reality there will be less demand than expected due to the higher prices: thus resulting in a profit loss.

However, the actual profits of the MNL are slightly higher with the heuristic than with the exact method. This is due to the frequencies on the Rotterdam-Duisburg service being slightly higher in general with the MNL than with the Mixed Logit. Since the actual shippers population has a frequency coefficient $\beta_f^{\mathrm{Inter}}$ higher than both the MNL and Mixed Logit (see Table~\ref{tab:coeffU}), the decision to propose additional frequency pays off.

So far, we notice that the heuristic is able to find solutions performing at least as well as the exact method. We can then make use of the fast solving time of the heuristic to increase the number of draws from 1000 to 5000. Table~\ref{tab:resheur5K} reports the solutions' statistics with a sample size of $R = 5000$.

\begin{table}[!t]
\centering
\caption{Solutions of stochastic models with predetermination heuristic (5000 draws)}
\label{tab:resheur5K}
\resizebox{0.9\textwidth}{!}{%
\begin{tabular}{llccc|ccc}
\multicolumn{1}{c}{} & \multicolumn{1}{c}{} & \multicolumn{3}{c}{MNL} & \multicolumn{3}{c}{Mixed Logit} \\ \hline
\multicolumn{1}{c}{} & \multicolumn{1}{c}{} & Min. & Average & Max. & Min. & Average & Max. \\ \hline
\multirow{4}{*}{\begin{tabular}[c]{@{}l@{}}Weekly\\ frequencies\\ {[}M8 vessels \\ (M11 vessels){]}\end{tabular}} & RTM-DUI & 0 (0) & 0 (10) & 0 (14) & 0 (0) & 0 (5) & 0 (13) \\
 & RTM-BON & 0 (0) & 0 (0) & 0 (0) & 0 (0) & 0 (0) & 0 (0) \\
 & DUI-BON & 0 (0) & 7 (2) & 12 (2) & 0 (1) & 10 (1) & 16 (3) \\
 & RTM-DUI-BON & 21 (0) & 22 (0) & 24 (0) & 19 (0) & 21 (0) & 24 (0) \\ \hline
\multirow{6}{*}{Prices {[}\euro{}{]}} & RTM-DUI & 183 & 198 & 222 & 149 & 163 & 176 \\
 & DUI-RTM & 174 & 191 & 203 & 138 & 155 & 177 \\
 & RTM-BON & 181 & 203 & 244 & 159 & 180 & 201 \\
 & BON-RTM & 165 & 190 & 221 & 150 & 172 & 196 \\
 & DUI-BON & 137 & 158 & 187 & 111 & 142 & 169 \\
 & BON-DUI & 141 & 160 & 187 & 111 & 139 & 166 \\ \hline
\multicolumn{2}{l}{Computation time {[}hours{]}} & 0.38 & 0.47 & 0.55 & 0.52 & 0.55 & 0.60 \\ \hline
\end{tabular}%
}
\end{table}

There is a direct proportionality between the computation time and the number of random draws, suggesting a linear computational complexity of the predetermination heuristic. This would allow to further increase the number of draws, but some memory issues can occur as the precomputation step requires to store a lot of values before the optimization. Nevertheless, moving to 5000 draws allows to improve the accuracy of the returned solutions, as the interval between the minimal and maximal values of the decision variables becomes tighter.

\begin{figure}[!b]
\centering
\includegraphics[width=0.8\textwidth]{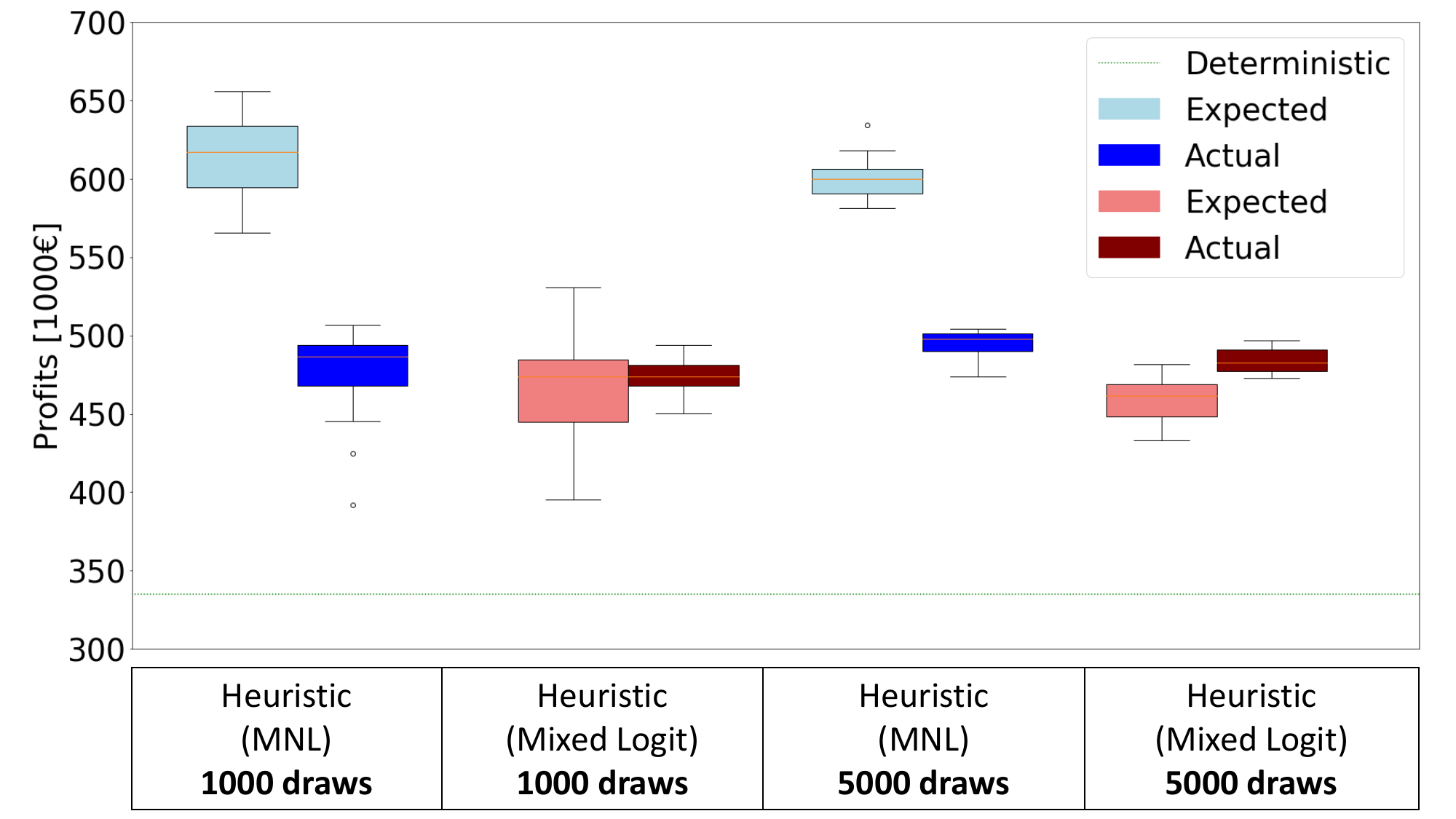}
\caption{Comparison of profits by the stochastic models using different number of draws with the heuristic.\label{fig:HeuristicDraws}}
\end{figure}

As depicted in Figure~\ref{fig:HeuristicDraws}, increasing the number of draws also leads to more accurate expected and actual profits since the confidence interval is tighter with 5000 draws than with 1000 draws. Stated differently, increasing the number of draws makes the heuristic less sensitive to changes of the sample, leading to more stable solutions. Additionally, the actual profits are enhanced when more draws are performed. There is still a significant drop between the expected and actual profits in the MNL case because the model underestimates the cost sensitivity of shippers compared to reality. For the Mixed Logit, there is a noticeable increase in the actual profits compared to the expected ones. Referring back to Table~\ref{tab:coeffU}, the average value of the cost coefficient $\bar{\beta}_c^{\mathrm{Inter}}$ is higher in absolute value for the Mixed Logit than in the actual population. This means that the CD-SNDP with Mixed Logit has the tendency to overestimate the cost senstivity of shippers compared to reality. Therefore, it will be more cautious in pricing but will in the end be able to attract more demand than anticipated, which results in increased profits.

\subsection{Key insights}\label{recommandations}

Several take-aways can be gathered from the results presented above. First, a cycle-based formulation (with multiple stops allowed) of the SND problem is more efficient in terms of asset usage as the operator can use consolidation opportunities. This results in both reduced costs and increased demand. The mathematical expression of services is less straightforward than with a path-based formulation due to the addition of service legs, but the improved results justify this effort.

Secondly, it is highly beneficial for the transport operator to include the information they have about the demand during the design of their services. The CD-SNDP results have shown that, even with a simple deterministic model, the solution of the SNDP problem is able to generate actual profits that are nearly three times higher than the benchmark. This is because the benchmark's assumption that shippers are purely cost-minimizers neglects other attributes that still play a role in the decision-making of shippers, such as the service frequencies. The utility functions also include the arbitrage between these attributes through the weighting coefficients.

Thirdly, making use of stochastic CD-SNDP exploits further the potential of the model. Indeed, perfect and complete information about the shippers is not available to the operator, so that their demand model will miss some aspects that play a role in the shippers' choices. These aspects can indirectly be accounted for by adding random error terms in the model. Including this uncertainty into the model enables gains of almost 50\% compared to the deterministic CD-SNDP. Therefore, the stochastic formulation of the CD-SNDP is one convenient way to account for imperfect information endogeneously to the model.

Finally, quantifying and incorporating the heterogeneous preferences of shippers allows for a more accurate estimation of the profits. Indeed, except for the stochastic CD-SNDP with Mixed Logit, all models presented above substantially overestimate the profits. This can lead to very bad surprises for the operator if they expect a given amount of profit in their budget, but end up realizing much less. On the other hand, the formulation with Mixed Logit expects profits in line with (or even lower than) the ones that are realized. Considering heterogeneity then allows to get a better prevision of the profits. 

\section{Conclusion}\label{conclu}

This work proposes a Service Network Design and Pricing problem that incorporates the mode choice behavior of shippers. Therefore, we develop a so-called Choice-Driven Service Network Design and Pricing problem that directly includes utility-based mode choice models into a bilevel optimization problem, which can then be reformulated as a single level linear problem. The random nature of utility-based models, such as the Multinomial Logit, allows to account for missing information about attributes playing a role in the mode choice. Opting for a Mixed Logit formulation further allows to consider the heterogeneous preferences of shippers, thus getting a more realistic representation of the shippers' population. Due to the randomness, the problem becomes stochastic, which makes it computationally expensive to solve with an exact method. To overcome this issue, we develop a predetermination heuristic that computes utilities prior to the optimization.

The results show that the heuristic is able to considerably reduce the computational time, while finding solutions of similar quality to the exact method. Regarding the proposed model itself, it is compared to a benchmark where shippers are assumed purely cost-minimizers. We show that the profits achieved by our model are substantially higher, even if the embedded mode choice model is simply deterministic. All in all, including more information about the shippers while designing and pricing the services suggests considerable gains for the transport operator. Even if the exact model or parameters are not known, it is still far better than not using the available information.

Now that the potential of our Choice-Driven Service Network Design and Pricing has been demonstrated on a small network, the immediate next step is to apply this methodology to a larger network. This will allow to assess the proposed model on a real-size instance and to test the performance of the predetermination heuristic on a larger problem. The search space could potentially be reduced by applying practical rules, such as a maximum number of stops per services.

Moreover, several assumptions made in this work deserve to be challenged. Firstly, in the mode choice models, the utilities of the vessel operator and of the competing IWT carrier are considered independent from each other. However, since they are both proposing inland waterway services, these two options are correlated with each other. This could also apply to a lesser extent to the Rail alternative, which also proposes scheduled intermodal services. Further work should consider this correlation between choice alternatives. Secondly, it is assumed that the utilities of the competing alternatives can be computed by the operator, implying that they have full information about their competitors. Some attributes can indeed be found, e.g., the frequency or travel times, but the price that the competitors are applying cannot be known perfectly: at best it can be estimated. The choice-driven model should then be developed further to account for this imperfect information. Thirdly, the competition is assumed exogenous and fixed meaning that they will not react to the operator's new services. But the competitors will also seek to improve their services and profits, even more so if they lose market share to the operator. These dynamics can be covered, for example, through an Agent-Based Model accounting for the reactions of the different parties involved.

Another dynamic aspect that can be included is about the pricing, particularly in the context of inland waterway transport. Indeed, the frequent low water levels on the Rhine reduce the capacity of the vessels, thus increasing the transportation costs per container. The operators then have to increase the price they charge to compensate for the losses. A time dimension could be included in the optimization model to deal with dynamic pricing. Finally, our formulation implies that a single price per OD pair is set for all shippers. However, revenue management techniques can be used to improve the performance of the proposed model. This will provide the operator with additional gains because they can tailor the prices offered to specific customers. A revenue management setting would also to develop the full potential of the formulation with Mixed Logit, as the prices can be adapted to the different cost sensitivities of shippers.

\ACKNOWLEDGMENT{%
This research is supported by the project “Novel inland waterway transport concepts for moving freight effectively (NOVIMOVE)”. This project has received funding from the European Union’s Horizon 2020 research and innovation programme under grant agreement No 858508.
}

%
%
%
\begin{APPENDICES}
\end{APPENDICES}


\bibliographystyle{informs2014trsc} 
\bibliography{SND.bib} 


\end{document}